\documentclass[11pt]{amsart}

\usepackage{amsmath}

\usepackage{amssymb}

\newcommand{\C}{{\mathbb C}}       
\newcommand{\R}{{\mathbb R}}       
\newcommand{\DD}{{\mathcal D}}
\newcommand{\HH}{{\mathcal H}}
\newcommand{\LL}{{\mathcal L}}
\newcommand{\D}{{\Delta}}

\newcommand{\RR}{{\mathcal R}}
\newcommand{\CH}{{\mathcal Ch}}
\newcommand{\EE}{{\mathcal E}}

\newcommand{\CC}{{\mathcal C}}
\newcommand{\CE}{{{\mathcal C}_{\varepsilon}}}
\newcommand{\CT}{{\wt{\CC}_{\varepsilon}}}

\newcommand{\diam}{{\rm diam}}
\newcommand{\dist}{{\rm dist}}
\newcommand{\ds}{\displaystyle }

\newcommand{\interior}[1]{{\stackrel{\mbox{\scriptsize$\circ$}}{#1}}}
\newcommand{\lra}{{\longrightarrow}}

\newcommand{\rf}[1]{{(\ref{#1})}}

\newcommand{\supp}{{\rm supp}}

\newcommand{\vphi}{{\varphi}}
\newcommand{\ve}{{\varepsilon}}

\newcommand{\vvv}{{\vspace{3mm}}}
\newcommand{\wt}[1]{{\widetilde{#1}}}
\newcommand{\wh}[1]{{\widehat{#1}}}

\newcommand{\noi}{\noindent}

\newtheorem{theorem}{Theorem}[section]
\newtheorem{lemma}[theorem]{Lemma}

\theoremstyle{definition}

\newtheorem{example}[theorem]{Example}

\theoremstyle{remark}
\newtheorem{rem}[theorem]{Remark}

\numberwithin{equation}{section}

\newcommand{\brem}{\begin{rem}}
\newcommand{\erem}{\end{rem}}

\newcommand{\bexam}{\begin{example}}
\newcommand{\eexam}{\end{example}}

\begin{document}

\title[Painlev\'e's problem and analytic
capacity]{Painlev\'e's problem and
the semiadditivity of analytic capacity}

\author[XAVIER TOLSA]{Xavier Tolsa}

\address{D\'epartement de Math\'ematique, Universit\'e de Paris-Sud,
91405 Orsay Cedex, France}

\curraddr{Departament de Matem\`atiques, Universitat Aut\`onoma de
Barcelona, 08193 Bellaterra (Barcelona), Spain}

\email{xtolsa@mat.uab.es}

\thanks{Supported by a Marie Curie Fellowship of the European Community
program Human Potential under contract HPMFCT-2000-00519.
Also partially supported by grants DGICYT BFM2000-0361 (Spain)
and 2000/SGR/00059 (Generalitat de Catalunya).}


\date{September, 2001.}

\begin{abstract}
Let $\gamma(E)$ be the analytic capacity of a compact set $E$ and let
$\gamma_+(E)$ be the capacity of $E$ originated by
Cauchy transforms of positive measures. In this paper we prove
that $\gamma(E)\approx\gamma_+(E)$ with estimates independent of $E$. As a
corollary, we characterize removable singularities for bounded analytic
functions in terms of curvature of measures, and
we deduce that $\gamma$ is semiadditive.
\end{abstract}

\maketitle




\section{Introduction}

The {\em analytic capacity} of a compact set $E\subset\C$ is defined as
$$\gamma(E) = \sup|f'(\infty)|,$$
where the supremum is taken over all analytic functions $f:\C\setminus E\lra
\C$ with $|f|\leq1$ on $\C\setminus E$, and $f'(\infty)=\lim_{z\to\infty}
z(f(z)-f(\infty))$. For a general set $F\subset \C$, we set
$\gamma(F)= \sup\{\gamma(E):E\subset F,\mbox{ E compact}\}.$

The notion of analytic capacity was first introduced by Ahlfors \cite{Ahlfors}
in the 1940's in order
to study the removability of singularities of bounded analytic
functions. A compact set $E\subset \C$ is said to be removable (for
bounded analytic functions) if for any open set $\Omega$ containing $E$, every
bounded function analytic on $\Omega\setminus E$ has an analytic
extension to $\Omega$. In \cite{Ahlfors} Ahlfors
showed that $E$ is removable if and only if $\gamma(E)=0$. However, this
result doesn't characterize removable singularities for
bounded analytic functions in a geometric way, since the definition of
$\gamma$ is purely analytic.

Analytic capacity was rediscovered by Vitushkin in the 1950's in connection
with problems of uniform approximation of analytic functions by rational
functions (see \cite{Vitushkin}, for example).
He showed that analytic capacity plays a central role in this type of
problems. This fact motivated a renewed interest in analytic capacity.
The main drawback of Vitushkin's
techniques arises from the fact that there is not a complete description
of analytic capacity in metric or geometric terms.

On the other hand, the {\em analytic capacity $\gamma_+$}
(or {\em capacity $\gamma_+$}) of a compact set $E$ is
$$\gamma_+(E) = \sup_{\mu} \mu(E),$$
where the supremum is taken over all positive Radon measures $\mu$ supported
on $E$ such that the Cauchy transform $\ds f= \frac1z *\mu$ is an
$L^\infty(\C)$ function with $\|f\|_\infty\leq1$. Since
$\Bigl(\ds\frac1z * \mu\Bigr)'(\infty)=\mu(E)$, we have
\begin{equation} \label{eqtriv}
\gamma_+(E) \leq \gamma(E).
\end{equation}
To the best of our knowledge, the capacity $\gamma_+$ was introduced by Murai
\cite[pp. 71-72]{Murai}. He showed that some estimates on $\gamma_+$
are related to the $L^2$ boundedness of the Cauchy transform.

In this paper we prove the converse of inequality \rf{eqtriv} (modulo a
multiplicative constant):

\begin{theorem} \label{semiadd}
There exists an absolute constant $A$ such that
$$\gamma(E) \leq A\gamma_+(E)$$
for any compact set $E$.
\end{theorem}

Therefore, we deduce $\gamma(E) \approx \gamma_+(E)$ (where $a\approx b$
means that there exists an absolute positive constant $C$ such that
$C^{-1}b \leq a \leq Cb$), which was a quite old question concerning
analytic capacity (see for example \cite{DO} or \cite[p.435]{Verdera1.5}).

To describe the consequences of Theorem \ref{semiadd} for Painlev\'e's problem
(that is, the problem of characterizing removable singularities for bounded
analytic functions in a geometric way) and for the semiadditivity of analytic
capacity, we need to introduce some additional notation and terminology.

Given a complex Radon measure $\nu$ on $\C$,
the {\em Cauchy transform} of $\nu$ is
$$\CC \nu(z) =\int \frac{1}{\xi-z}\, d\nu(\xi).$$
This definition does not make sense, in general, for
$z\in\supp(\nu)$, although one can easily see that the integral above is
convergent at a.e. $z\in\C$ (with respect to Lebesgue measure).
This is the reason why one considers the {\em truncated
Cauchy transform} of $\nu$, which is defined as
$$\CE \nu(z) =\int_{|\xi-z|>\ve} {\frac1{\xi-z}}\, d\nu(\xi),$$
for any $\ve>0$ and $z\in\C$.
Given a $\mu$-measurable function $f$ on $\C$ (where $\mu$ is some fixed
positive Radon measure on $\C$), we write
$\CC f \equiv \CC (f \,d\mu)$
and
$\CE f \equiv \CE (f \,d\mu)$
for any $\ve>0$.
It is said that the Cauchy transform is bounded on $L^2(\mu)$ if the operators
$\CE$ are bounded on $L^2(\mu)$ uniformly on $\ve>0$.

A positive Radon measure $\mu$ is said to have linear growth if there exists
some constant $C$ such that $\mu(B(x,r))\leq Cr$ for all $x\in\C$, $r>0$.

Given three pairwise different points $x,y,z\in\C$, their {\em Menger
curvature} is
$$c(x,y,z) = \frac{1}{R(x,y,z)},$$
where $R(x,y,z)$ is the radius of the circumference passing through $x,y,z$
(with $R(x,y,z)=\infty$, $c(x,y,z)=0$ if $x,y,z$ lie on a same line). If two
among these points coincide, we let $c(x,y,z)=0$.
For a positive Radon measure $\mu$, we set
$$c^2_\mu(x) = \int\!\!\int c(x,y,z)^2\, d\mu(y) d\mu(z),$$
and we define the {\em curvature of $\mu$} as
\begin{equation} \label{defcurv}
c^2(\mu) = \int c^2_\mu(x)\, d\mu(x)
= \int\!\!\int\!\!\int c(x,y,z)^2\, d\mu(x) d\mu(y) d\mu(z).
\end{equation}
On the one hand, the notion of curvature of a measure,
first introduced by Melnikov \cite{Melnikov} when he was studying a discrete
version of analytic capacity, is connected to the Cauchy
transform. This relationship
stems from the following remarkable identity found by Melnikov and Verdera
\cite{MV} (assuming that $\mu$ has linear growth):
\begin{equation}\label{idmv}
\|\CE\mu\|_{L^2(\mu)}^2 = \frac{1}{6} c^2_\ve(\mu) +
O(\mu(\C)),
\end{equation}
where $c^2_\ve(\mu)$ is an $\ve$-truncated version of
$c^2(\mu)$ (defined as in the right hand side of \rf{defcurv}, but with the
integrals over $\{x,y,z\in\C: |x-y|,|y-z|,|x-z|>\ve\}$). On the other
hand, the curvature of a measure encodes metric and geometric
information from the support of the measure and is related to
rectifiability (see \cite{Leger}). In fact, there is a
close relationship between $c^2(\mu)$ and the coefficients $\beta$ which
appear in Jones' traveling salesman result \cite{Jones}.

Using the identity \rf{idmv},
it has been shown in \cite{Tolsa1} that the capacity $\gamma_+$ has a
rather precise description in terms of curvature of measures
(see \rf{cap++2} and \rf{cap++4}). As a consequence,
from Theorem \ref{semiadd} we get a characterization of analytic capacity
with a definite metric-geometric flavour. In particular,
in connection with Painlev\'e's problem
we obtain the following result, previously conjectured by Melnikov
(see \cite{Davidsurvey} or \cite{Mattila3}).

\begin{theorem}\label{painleve1}
A compact set $E\subset\C$ is non removable for bounded analytic functions
if and only if it supports a positive Radon measure with linear growth
and finite curvature.
\end{theorem}

It is easy to check that this result follows from the
comparability between $\gamma$ and $\gamma_+$. In fact,
it can be considered as a qualitative version of Theorem
\ref{semiadd}.

From Theorem \ref{semiadd} and \cite[Corollary 4]{Tolsa2}
we also deduce the following result,
which in a sense can be considered as the dual of Theorem \ref{painleve1}.

\begin{theorem} \label{painleve2}
A compact set $E\subset\C$ is removable for bounded analytic functions
if and only if there exists a finite positive Radon measure $\mu$ on $\C$ such
that for all $x\in E$ either $\Theta^*_\mu(x)=\infty$ or $c^2_\mu(x) =\infty$.
\end{theorem}

In this theorem, $\Theta^*_\mu(x)$ stands for the upper linear density of
$\mu$ at $x$, i.e. $\Theta^*_\mu(x)= \limsup_{r\to0} \mu(B(x,r))\,r^{-1}$.

Theorem \ref{semiadd} has another important corollary. Up to now, it was not
known if analytic capacity is semiadditive, that is, if there exists an
absolute constant $C$ such that
\begin{equation} \label{eqsem}
\gamma(E\cup F) \leq C(\gamma(E) + \gamma(F)).
\end{equation}
This question was raised by Vitushkin in the early 1960's (see
\cite{Vitushkin} and \cite{VM})
and was known to be true only in some particular cases
(see \cite{Melnikov1} and \cite{Suita} for example, and
\cite{Davie} and \cite{DO} for some related results).
On the other hand, a positive answer to the semiadditivity problem would
have interesting applications to rational approximation (see \cite{Verdera1.5}
and \cite{Vitushkin}). Theorem \ref{semiadd} implies that, indeed, analytic
capacity is semiadditive because $\gamma_+$ is semiadditive
\cite{Tolsa1}. In fact, the following stronger result holds.

\begin{theorem} \label{corosemiadd}
Let $E\subset\C$ be compact. Let $E_i$, $i\geq1$, be Borel sets such
that $E=\bigcup_{i=1}^\infty E_i$. Then,
$$\gamma(E) \leq C\sum_{i=1}^\infty \gamma(E_i),$$
where $C$ is an absolute constant.
\end{theorem}

Several results dealing with analytic capacity have been obtained
recently. Curvature of measures plays an essential role in most of them.
G. David proved in \cite{David} that a compact set $E$ with finite
length, i.e. with
$\HH^1(E)<\infty$ (where $\HH^s$ stands for the $s$-dimensional Hausdorff
measure), has vanishing analytic capacity if and only if it is purely
unrectifiable, that is, if $\HH^1(E\cap \Gamma)=0$ for all rectifiable
curves $\Gamma$.
This result had been known as Vitushkin's conjecture for a long time.
Let us also mention that in \cite{MMV} the same result had been proved
previously under an additional regularity assumption on the set $E$.

David's theorem is a very remarkable result.
However, it only applies to sets
with finite length. Indeed, Mattila \cite{Mattila1} showed that the natural
generalization of Vitushkin's conjecture to sets with non $\sigma$-finite
length does not hold (see also \cite{JM}).

After David's solution of Vitushkin's conjecture, Nazarov, Treil and Volberg
\cite{NTV} proved a $T(b)$ theorem useful for dealing with analytic
capacity. Their theorem also solves (the last step of) Vitushkin's
conjecture. Moreover, they obtained some quantitative results which imply
the following estimate:
\begin{equation} \label{nnntv}
\gamma_+(E) \geq C^{-1} \gamma(E)\,\biggl(1+
\left(\frac{\diam(E)}{\gamma(E)}\right)^2
\left(\frac{\HH^1(E)}{\gamma(E)}\right)^{38}\biggr)^{-1/2}.
\end{equation}
Notice that if $\HH^1(E)=\infty$, then the right hand side equals 0, and so
this inequality is not useful in this case.

On the other hand, some problems related to the capacity $\gamma_+$
have been studied recently. Some density estimates for $\gamma_+$ (among
other results) have been obtained in \cite{MP}, while in \cite{Tolsa2} it has
been shown that $\gamma_+$ verifies some properties which usually hold
for other capacities generated by positive potentials and energies, such as
Riesz capacities. Now all these results apply automatically to analytic
capacity, by Theorem \ref{semiadd}. See also \cite{MP1} and \cite{VMP}
for other questions related to $\gamma_+$.

Let us mention some additional consequences of Theorem \ref{semiadd}.
Up to now it was not even known if the class of sets with vanishing analytic
capacity was invariant under affine maps such as $x+iy\mapsto x+i2y$,
$x,y\in\R$ (this question was raised by O'Farrell, as far as we know).
However, this is true for $\gamma_+$ (and so for
$\gamma$), because its characterization
in terms of curvature of measures. Indeed, it is quite easy to check that
the class of sets with vanishing capacity $\gamma_+$ is
invariant under $\CC^{1+\ve}$ diffeomorphisms (see \cite{Tolsatesi}, for
example). The analogous fact for
$\CC^1$ or bilipschitz diffeomorphisms is an open problem.

Also, our results imply that David's theorem can be extended to sets
with $\sigma$-finite length. That is, if $E$ has $\sigma$-finite
length, then $\gamma(E) =0$ if and only if $E$ is purely
unrectifiable. This fact, which also remained unsolved,
follows directly either from Theorem \ref{semiadd} or \ref{corosemiadd}.

The proof of Theorem \ref{semiadd} in this paper is inspired by the recent
arguments of \cite{MTV},
where it is shown that $\gamma$ is comparable to $\gamma_+$ for a big class
of Cantor type sets. One essential new idea from \cite{MTV} is
the ``induction on scales'' technique, which can be also adapted to
general sets, as we shall see. Let us also remark that another important
ingredient of the proof of Theorem \ref{semiadd} is the $T(b)$ theorem of
\cite{NTV}.

Theorems \ref{painleve1} and \ref{painleve2} follow easily from Theorem
\ref{semiadd} and known results about $\gamma_+$. Also, to prove Theorem
\ref{corosemiadd}, one only has to use Theorem \ref{semiadd} and the fact
that $\gamma_+$ is countably semiadditive. This has been shown in
\cite{Tolsa1} under the additional assumption that the sets $E_i$ in Theorem
\ref{corosemiadd} are compact. With some minor modifications, the proof
in \cite{Tolsa1} is also valid if the sets $E_i$ are Borel.
For the sake of completeness, the detailed arguments are shown in Remark
\ref{compness}.

The plan of the paper is the following. In Section \ref{prelim} we introduce
some notation and recall some preliminary results.
In Section \ref{sketch}, for the reader's convenience, we sketch the
ideas involved in the proof of Theorem \ref{semiadd}.
In Section \ref{sec3} we
prove a preliminary lemma which will be necessary for Theorem \ref{semiadd}.
The rest of the paper is devoted to
the proof of this theorem, which we have split into three parts.
The first one corresponds to the First Main Lemma \ref{main1}, which is stated
in Section
\ref{sec5} and proved in Sections \ref{sec6}--\ref{excep}. The second one
is the Second Main Lemma \ref{main2}, stated in Section \ref{sec8} and proved
in Sections
\ref{secT}--\ref{secsup}. The last part of the proof of Theorem \ref{semiadd}
is in Section \ref{sec11} and consists of the induction argument on scales.


\section{Notation and background} \label{prelim}

We denote by $\Sigma(E)$ the set of all positive Radon measures $\mu$
supported on $E\subset\C$ such that $\mu(B(x,r))\leq r$ for all $x\in E$,
$r>0$.

As mentioned in the Introduction, curvature of measures was
introduced by Melnikov in
\cite{Melnikov}. In this paper he proved the following inequality:
\begin{equation} \label{Mark}
\gamma(E)\geq C \sup_{\mu\in\Sigma(E)}
\frac{\mu(E)^2}{\mu(E) + c^2(\mu)},
\end{equation}
where $C>0$ is some absolute constant.
In \cite{Tolsa1} it was shown that inequality \rf{Mark} also holds
if one replaces $\gamma(E)$ by $\gamma_+(E)$ on the left hand side,
and then one obtains
\begin{equation}
\hspace{8mm}\gamma_+(E)
 \approx  \sup_{\mu\in\Sigma(E)} \frac{\mu(E)^2}{\mu(E) +
 c^2(\mu)}.
\label{cap++2}
\end{equation}

Let $M$ be the maximal radial Hardy-Littlewood operator:
$$M\mu(x) = \sup_{r>0} \frac{\mu(B(x,r))}{r}$$
[if $\mu$ were a complex measure, we would replace $\mu(B(x,r))$ by
$|\mu|(B(x,r))$],
and let $c_\mu(x) = \left(c_\mu^2(x)\right)^{1/2}$.
The following potential was introduced by Verdera in \cite{Verdera2}:
\begin{equation} \label{potencial}
U_\mu(x) := M\mu(x) + c_\mu(x),
\end{equation}
It turns out that $\gamma_+$ can also be characterized in terms of this
potential (see \cite{Tolsa2}, and also \cite{Verdera2} for a related result):
\begin{equation} \label{cap++4}
\hspace{8mm}\gamma_+(E)
\approx  \sup \{\mu(E):\,\supp(\mu)\subset E,\,
U_\mu(x)\leq 1\,\forall x\in E \}.
\end{equation}
Let us also mention that the potential $U_\mu$ will be very important for the
proof of Theorem \ref{semiadd}.

\vvv
\begin{rem}  \label{compness}
Let us see that Theorem \ref{corosemiadd} follows easily from Theorem
\ref{semiadd} and the characterization \rf{cap++4} of $\gamma_+$.
Indeed, if $E\subset\C$ is compact and $E=\bigcup_{i=1}^\infty E_i$, where
$E_i$, $i\geq1$, are Borel sets, then we take a Radon measure $\mu$
such that $\gamma_+(E) \approx\mu(E)$ and $U_\mu(x)\leq 1$ for all $x\in
E$. For each $i\geq1$, let $F_i\subset E_i$ be compact and such that
$\mu(F_i) \geq \mu(E_i)/2$. Since $U_{\mu|F_i}(x)\leq 1$ for all
$x\in F_i$, we deduce $\gamma_+(F_i) \geq C^{-1}\mu(F_i)$. Then, from Theorem
\ref{semiadd} we get
\begin{eqnarray*}
\gamma(E) & \approx & \gamma_+(E) \,\approx\, \mu(E) \,\leq\, C \sum_i
\mu(F_i)\\
& \leq & C \sum_i \gamma_+(F_i)
\,\approx\, C \sum_i \gamma(F_i) \,\leq\, C \sum_i \gamma(E_i).
\end{eqnarray*}
\end{rem}

\vvv
Let us recall the definition of the maximal Cauchy transform of a complex
measure $\nu$:
$$\CC_*\nu(x) = \sup_{\ve>0} |\CE\nu(x)|.$$
Let $\psi$ be a $\CC^\infty$ radial function supported on
$B(0,1)$, with $0\leq\psi\leq2$, $\|\nabla\psi\|_\infty\leq10$,
and $\int\psi d\LL^2=1$ (where $\LL^2$ stands for the Lebesgue
measure). We denote $\psi_\ve(x) = \ve^{-2}\psi(x/\ve)$. The
regularized operators $\CT$ are defined as
$$\CT \nu :=\psi_\ve *\CC\nu = \psi_\ve * \frac1z *\nu.$$
Let $r_\ve = \psi_\ve*\ds\frac1z$. It is easily seen that
$r_\ve(z)=1/z$ if $|z|>\ve$, $\|r_\ve\|_\infty\leq C/\ve$, and
$|\nabla r_\ve(z)|\leq C|z|^{-2}$. Further, since $r_\ve$ is a
uniformly continuous kernel, $\CT\nu$ is a continuous function on
$\C$. Notice also that if $|\CC\nu|\leq B$ a. e. with respect to
Lebesgue measure, then $|\CT(\nu)(z)|\leq B$ for all $z\in\C$.

Moreover, we have
\begin{equation} \label{w1}
|\CT\nu(x)-\CE\nu(x)| = \left| \int_{|y-x|<\leq\ve} r_\ve(y-x)
d\nu(y)\right| \leq C\, M \nu(x).
\end{equation}

By a square $Q$ we mean a closed square with sides parallel to the axes.

The constant $A$ in Theorem \ref{semiadd} will be fixed at the end of
the proof. Throughout all the paper,
the letter $C$ will stand for an absolute constant
that may change at different occurrences. Constants with subscripts, such as
$C_1$, will retain its value, in general. On the other hand, the
constants $C,\,C_1,\ldots$ do {\bf not} depend on $A$.


\section{Outline of the arguments for the proof of Theorem
  \ref{semiadd}} \label{sketch}

In this section we will sketch the
arguments involved in the proof of Theorem \ref{semiadd}.

In the rest of the paper, unless stated otherwise, {\bf we will assume that
$E$ is a finite union of compact disjoint segments.}
We will prove Theorem \ref{semiadd} for this type of sets.
The general case follows from this particular instance by a
discretization argument, such as in \cite[Lemma 1]{Melnikov}.
Moreover, we will assume that the segments make an angle of $\pi/4$, say, with
the $x$ axis. In this way, the intersection of $E$ with any parallel line to
one the coordinate axes will always have zero length. This fact will avoid
some technical problems.

To prove Theorem \ref{semiadd} we want to apply some kind of $T(b)$
theorem, as David in \cite{David} for the proof of Vitushkin's conjecture.
Because of the definition of analytic capacity, there exists a
complex Radon measure $\nu_0$ supported on $E$ such that
\begin{eqnarray}
\|\CC\nu_0\|_\infty & \leq & 1 \label{cc1},\\
|\nu_0(E)| & = & \gamma(E)  \label{cc2},\\
d\nu_0 & = & b_0\,d\HH^1|E,\quad \mbox{with $\|b_0\|_\infty\leq 1$}.
 \label{cc3}
\end{eqnarray}
We would like to show that there exists some Radon measure $\mu$
supported on $E$ with $\mu\in\Sigma(E)$, $\mu(E) \approx \gamma(E)$,
and such that the Cauchy transform is bounded on $L^2(\mu)$ with
absolute constants. Then, using \rf{cap++2} for example, we would get
$$\gamma_+(E) \geq C^{-1}\mu(E)\geq C^{-1}\gamma(E),$$
and we would be done.

However, by a more or less direct application of a $T(b)$ theorem we
cannot expect to
prove that the Cauchy transform is bounded with respect to such a
measure $\mu$ with absolute constants. Let us explain the reasons in
some detail. Suppose for example that there exists some function $b$ such
that $d\nu_0 = b\,d\mu$ and we use the information
about $\nu_0$ given by \rf{cc1}, \rf{cc2} and \rf{cc3}. From \rf{cc1}
and \rf{cc2} we derive
\begin{equation}
\|\CC(b\,d\mu)\|_\infty \leq 1 \label{cc1'}
\end{equation}
and
\begin{equation}
\left|\int b\,d\mu\right| \approx \mu(E).\label{cc2'}
\end{equation}
The estimate \rf{cc1'} is very good for our purposes. In fact, most classical
$T(b)$-type theorems require only the $BMO(\mu)$ norm of $b$ to be bounded,
which is a weaker assumption. On the other hand, \rf{cc2'} is a global
paraaccretivity condition, and with some technical difficulties
(which may involve some kind of stopping time argument, like
in \cite{Christ}, \cite{David} or \cite{NTV}), one can
hope to be able to prove that the local paraaccretivity condition
$$\left|\int_Q b\,d\mu\right| \approx \mu(Q\cap E)$$
holds for many squares $Q$.

Our problems arise from \rf{cc3}. Notice that \rf{cc3} implies
that $|\nu_0|(E) \leq \HH^1(E)$, where $|\nu_0|$ stands for the variation of
$\nu_0$. This is a very bad estimate since we don't have any control on
$\HH^1(E)$ (we only know $\HH^1(E)<\infty$ because our assumptions on $E$).
However, as far as we know, all $T(b)$-type theorems require the estimate
\begin{equation}  \label{cc3'}
|\nu_0|(E)\leq C\mu(E)
\end{equation}
(and often stronger assumptions involving the $L^\infty$ norm of $b$).
So by a direct application of a $T(b)$-type theorem we will obtain
bad results when $\gamma(E) \ll \HH^1(E)$, and at most we will get estimates
which involve the ratio $\HH^1(E)/\gamma(E)$, such as \rf{nnntv}.

To prove Theorem \ref{semiadd}, we need to work with a measure ``better''
than $\nu_0$, which we call $\nu$. This new measure will be a suitable
modification of $\nu_0$ with the required estimate for its variation.
To construct $\nu$, we operate as in \cite{MTV}. We consider a set $F$
containing $E$ made up of a finite disjoint union of squares:
$F =\bigcup_{i\in I} Q_i$. One should think that the squares $Q_i$
approximate $E$ at some ``intermediate scale''. For example, in the
case of the usual $1/4$ planar Cantor set of generation $n$
studied in \cite{MTV}, the squares $Q_i$ are the squares of
generation $n/2$. For each square $Q_i$, we take a complex measure
$\nu_i$ supported on $Q_i$ such that $\nu_i(Q_i)= \nu_0(Q_i)$ and
 $|\nu_i|(Q_i)=|\nu_i(Q_i)|$ (that is, $\nu_i$ will be a constant
 multiple of a positive measure). Then we set $\nu=\sum_i\nu_i$.
So, if the squares $Q_i$ are big enough,
the variation $|\nu|$ will be sufficiently small. On the
 other hand, the squares $Q_i$ cannot be too big, because we will need
\begin{equation} \label{derf}
\gamma_+(F)\leq C\gamma_+(E).
\end{equation}

In this way, we will have constructed a complex measure $\nu$
supported on $F$ satisfying
\begin{equation}\label{derf2}
|\nu|(F) \approx |\nu(F)| = \gamma(E).
\end{equation}
Taking a suitable measure $\mu$ such that $\supp(\mu)\supset \supp(\nu)$ and
$\mu(F) \approx \gamma(E)$, we will be ready for the application of a
$T(b)$ theorem. Indeed, notice that \rf{derf2} implies that $\nu$
satisfies a global paraaccretivity condition and that also the variation
$|\nu|$ is controlled. On the other hand, if we have been careful
enough, we will have also some useful estimates on $|\CC\nu|$, since $\nu$
is an approximation of $\nu_0$ (in fact, when defining $\nu$ in
the paragraph above, the measures $\nu_i$ have to be constructed in a
smoother way).
Using the $T(b)$ theorem of \cite{NTV}, we will deduce
$$\gamma_+(F) \geq C^{-1}\mu(E),$$
and so, $\gamma_+(E)\geq C^{-1}\gamma(E)$, by \rf{derf}, and we will be done.

Several difficulties arise in the implementation of the arguments
above. In order to obtain the right estimates on the measures
$\nu$ and $\mu$ we will need to assume that
$\gamma(E\cap Q_i) \approx \gamma_+(E\cap Q_i)$ for each
square $Q_i$. For this reason, we will use
an induction argument involving the sizes of the squares $Q_i$, as in
\cite{MTV}. Further, the choice of the right squares $Q_i$ which approximate
$E$ at an intermediate scale is more complicated than in \cite{MTV}.
A careful examination of the arguments in \cite{MTV} shows the following.
Let $\sigma$ be an
extremal measure for the right hand side of \rf{cap++2}, and so for
$\gamma_+$ in a sense (now $E$ is some precise planar Cantor set). It is not
difficult to check that $U_\sigma(x) \approx1$ for all $x\in E$
(remember \rf{potencial}). Moreover, one can also check that the
corresponding squares $Q_i$ satisfy
\begin{equation} \label{derf10}
U_{\sigma|2Q_i}(x)\approx U_{\sigma|\C\setminus 2Q_i}(x)\approx 1
\quad \mbox{ for all $x\in Q_i$}.
\end{equation}
In the general situation of $E$ given a finite union of disjoint compact
segments,
the choice of the squares $Q_i$ will be also determined by the potential
$U_\sigma$, where now $\sigma$ is corresponding maximal measure for
the right hand side of \rf{cap++2}. We will not ask the squares $Q_i$ to
satisfy \rf{derf10}. Instead we will use a variant of this idea.

Let us mention that the First Main Lemma \ref{main1} below deals with
the construction of the measures $\nu$ and $\mu$ and the estimates
involved in this construction. The Second Main Lemma \ref{main2}
is devoted to the application of a suitable $T(b)$ theorem.


\section{A preliminary lemma}  \label{sec3}

In next lemma we show a property of the capacity $\gamma_+$ and its associated
potential which will play an important role in the choice of the squares $Q_i$
mentioned in the preceding section.

\begin{lemma}  \label{extremal}
There exists a measure $\sigma\in \Sigma(E)$ such that
$\sigma(E)\approx\gamma_+(E)$ and $U_\sigma(x) \geq \alpha$ for all $x\in E$,
where $\alpha>0$ is an absolute constant.
\end{lemma}

Let us remark that a similar result has been proved in
\cite[Theorem 3.3]{Tolsa2}, but without the assumption
$\sigma\in\Sigma(E)$.

\begin{proof}
We will see first that there exists a Radon measure $\sigma\in\Sigma(E)$ such
that the supremum on the right hand side of \rf{cap++2} is attained by
$\sigma$. That is,
$$g(E) := \sup_{\mu\in\Sigma(E)} \frac{\mu(E)^2}{\mu(E) + c^2(\mu)}
 = \frac{\sigma(E)^2}{\sigma(E) + c^2(\sigma)}.$$
This measure will fulfill the required properties.

It is easily seen that any measure
$\mu\in\Sigma(E)$ can be written as
$d\mu = f\,d\HH^1|E$, with $\|f\|_{L^\infty(\HH^1|E)}\leq1$,
by the Radon-Nikodym theorem. Take a sequence of functions $\{f_n\}_n$, with
$\|f_n\|_{L^\infty(\HH^1|E)}\leq1$, converging weakly in $L^\infty(\mu)$ to
some function $f\in L^\infty(\mu)$ and such that
$$\lim_{n\to\infty}\frac{\mu_n(E)^2}{\mu_n(E) + c^2(\mu_n)} = g(E),$$
with $d\mu_n=f_n\,d\HH^1|E$, $\mu_n\in\Sigma(E)$.
Consider the measure $d\sigma=f\,d\HH^1|E$. Because of the weak convergence,
$\mu_n(E)\to\sigma(E)$ as $n\to\infty$, and moreover $\sigma\in\Sigma(E)$.
On the other hand, it is an easy exercise to check that $c^2(\sigma) \leq
\liminf_{n\to\infty} c^2(\mu_n)$. So we get
$$g(E) =  \frac{\sigma(E)^2}{\sigma(E) + c^2(\sigma)}.$$

Let us see that $\sigma(E)\approx\gamma_+(E)$. Since $\sigma$ is maximal and
$\sigma/2$ is also in $\Sigma(E)$, we have
$$\frac{\sigma(E)^2}{\sigma(E) + c^2(\sigma)}\geq
\frac{\sigma(E)^2/4}{\sigma(E)/2 + c^2(\sigma)/8}.$$
Therefore,
$$\frac{\sigma(E)}{2}+ \frac{c^2(\sigma)}{8} \geq
\frac{\sigma(E)}{4}+ \frac{c^2(\sigma)}{4}.$$
That is, $c^2(\sigma) \leq 2\sigma(E)$. Thus,
$$\gamma_+(E) \approx g(E) \approx \sigma(E).$$

It remains to show that there exists some $\alpha>0$ such that
$U_\sigma(x)
\geq \alpha$ for all $x\in E$.
Suppose that $M\sigma(x) \leq 1/1000$ for some $x\in E$, and let $B:=B(x,R)$ be
some fixed ball.
We will prove the following:

\vspace{4mm}

\noi {\bf Claim.} {\em If $R>0$ is small enough, then there exists some set
 $A\subset B(x,R)\cap E$, with $\HH^1(A)>0$, such that the measure
$\sigma_\lambda :=\sigma +
\lambda\HH^1|A$ belongs to $\Sigma(E)$ for $0\leq\lambda\leq 1/100$.}

\begin{proof}[Proof of the claim]
Since $E$ is made up of a finite number of disjoint compact segments, we may
assume that $R>0$ is so small that $\HH^1(B(y,r)\cap E)\leq2r$ for all
$y\in B$, $0<r\leq 4R$, and also that $\HH^1(B(x,R)\cap E)\geq R$.
These assumptions imply that for any subset $A\subset B$ we have
$$
\HH^1(A\cap B(y,r))\leq \HH^1(E\cap B(y,r) \cap B) \leq 2r \qquad \mbox{for
all $y\in B$, $r>0$}.
$$
Thus, $\HH^1(A\cap B(y,r))\leq 4r$ for all $y\in\C$ and so
\begin{equation} \label{mh1}
M(\HH^1|A)(y) \leq 4 \qquad \mbox{for all $y\in \C$.}
\end{equation}

We define $A$ as
$$A = \{y\in B: M\sigma(y) \leq 1/4\}.$$
Let us check that $\HH^1(A)>0$. Notice that
\begin{equation} \label{m2B}
\sigma(2B) \leq 2R\,M\sigma(x) \leq\frac{2R}{1000} \leq \frac{1}{500}\,\HH^1
(B\cap E).
\end{equation}
Let $D=B\setminus A$. If $y\in D$, then $M\sigma(y) > 1/4$.
If $r>0$ is such that $\sigma(B(y,r))/r >1/4$, then $r\leq R/10$. Otherwise,
$B(y,r)\subset B(x,12r)$ and so
$$\frac{\sigma(B(y,r))}{r} \leq \frac{\sigma(B(x,12r))}{r} \leq 12M\sigma(x)
\leq \frac{12}{1000}.$$
Therefore,
$$D\subset \{y\in B: M(\sigma|2B)(y) > 1/4\}.$$
For each $y\in D$,
take $r_y$ with $0<r_y\leq R/10$ such that
$$\frac{\sigma(B(y,r_y))}{r_y} >\frac14.$$
By Vitali's $5r$-covering Theorem there are some disjoint balls
$B(y_i,r_{y_i})$ such that $D\subset \bigcup_i B(y_i,5r_{y_i})$.
Since we must have $r_{y_i}\leq R/10$, we get
$\HH^1(B(y_i,5r_{y_i})\cap E) \leq 15r_{y_i}$. Then, by \rf{m2B} we deduce
\begin{eqnarray*}
\HH^1(D\cap E) & \leq & \sum_i \HH^1(B(y_i,5r_{y_i})\cap E) \\
& \leq & \sum_i 15r_{y_i} \leq 60 \sum_i \sigma(B(y_i,r_{y_i})) \\
& \leq & 60\sigma(2B) \leq \frac{60}{500} \HH^1(B\cap E).
\end{eqnarray*}
Thus, $\HH^1(A)>0.$

Now we have to show that $M\sigma_\lambda(y) \leq 1$ for all $y\in E$.
If $y\in A$, then $M\sigma(y) \leq1/4$, and then by \rf{mh1} we have
$$M\sigma_\lambda(y) \leq \frac14 + \lambda M(\HH^1|A)(y) \leq \frac14 +
\frac{4}{100} <1.$$
If $y\not\in A$ and $B(y,r)\cap A=\varnothing$, then we obviously have
$$\frac{\sigma_\lambda(B(y,r))}r = \frac{\sigma(B(y,r))}r \leq 1.$$
Suppose $y\not\in A$ and $B(y,r)\cap A\neq\varnothing$. Let $z\in
B(y,r)\cap A$. Then,
$$\frac{\sigma(B(y,r))}r \leq \frac{\sigma(B(z,2r))}r \leq 2M\sigma(z) \leq
\frac12.$$
Thus,
$$\frac{\sigma_\lambda(B(y,r))}r \leq \frac12 + \lambda M(\HH^1|A)(y)
\leq \frac12 + \frac{4}{100} <1.$$
So we always have $M\sigma_\lambda(y) \leq 1$.
\end{proof}

\vspace{3mm}
Let us continue the proof of Lemma \ref{extremal} and
let us see that $U_\sigma(x) \geq \alpha$. Consider the function
$$\vphi(\lambda) = \frac{\sigma_\lambda(E)^2}{\sigma_\lambda(E) +
c^2(\sigma_\lambda)}.$$
Since $\sigma$ is a maximal measure for $g(E)$ and
$\sigma_\lambda\in\Sigma(E)$ for some $\lambda>0$, we must have
$\vphi'(0)\leq 0$.
Observe that
\begin{eqnarray*}
\vphi(\lambda) & = &
\bigl[\sigma(E) + \lambda\HH^1(A)\bigr]^2\\
& & \times \Bigl[\sigma(E) +
\lambda\HH^1(A) + c^2(\sigma) + 3\lambda\, c^2(\HH^1|A,\sigma,\sigma) \\
&& \mbox{}+ 3\lambda^2 c^2(\sigma,\HH^1|A,\HH^1|A) +
\lambda^3 c^2(\HH^1|A) \Bigr]^{-1}.
\end{eqnarray*}
So,
$$\vphi'(0) = \frac{2 \sigma(E)\HH^1(A)(\sigma(E)+c^2(\sigma)) - \sigma(E)^2\,
(\HH^1(A) + 3\,c^2(\HH^1|A,\sigma,\sigma))}{(\sigma(E) + c^2(\sigma))^2}.$$
Therefore, $\vphi'(0)\leq 0$ if and only if
$$2\HH^1(A)(\sigma(E)+c^2(\sigma)) \leq \sigma(E)\,
(\HH^1(A)+ 3\,c^2(\HH^1|A,\sigma,\sigma)).$$
That is,
$$\frac{\sigma(E) + 2c^2(\sigma)}{\sigma(E)} \leq \frac{3\,
c^2(\HH^1|A,\sigma,\sigma)}{\HH^1(A)}.$$
Therefore,
$c^2(\HH^1|A,\sigma,\sigma)/\HH^1(A)\geq \frac{1}{3}$.
So there exists some $x_0\in A$ such that
\begin{equation} \label{wed11}
c^2(x_0,\sigma,\sigma)\geq \frac{1}{3}.
\end{equation}

We write
\begin{eqnarray} \label{wed12}
c^2(x_0,\sigma,\sigma) & = & c^2(x_0,\sigma|2B,\sigma|2B) +
c^2(x_0,\sigma|2B,\sigma|\C\setminus2B) \nonumber \\
&& \mbox{} + c^2(x_0,\sigma|\C\setminus2B,\sigma|\C\setminus2B).
\end{eqnarray}
If $R$ is chosen small enough, then $B\cap E$ coincides with a segment
and so $c^2(x_0,\sigma|2B,\sigma|2B) =0$. On the other hand,
\begin{eqnarray*}
c^2(x_0,\sigma|2B,\sigma|\C\setminus2B)  & \leq &
C\int_{y\in 2B}\int_{z\in\C\setminus 2B} \frac{1}{|x-z|^2}\,d\sigma(y)\,
d\sigma(z) \\
& \leq & C_2 M\sigma(x)^2.
\end{eqnarray*}
Thus, if $M\sigma(x)^2\leq 1/(6C_2)$, then  by \rf{wed11} and \rf{wed12} we
obtain
$$
c^2_{\sigma|\C\setminus 2B}(x_0) =
c^2(x_0,\sigma|\C\setminus2B,\sigma|\C\setminus2B) \geq
\frac{1}{3} - \frac16 = \frac16.
$$
Also, it is easily checked that
\begin{equation}\label{xx1}
|c_{\sigma|\C\setminus 2B}(x) - c_{\sigma|\C\setminus 2B}(x_0)|
\leq C_3M\sigma(x).
\end{equation}
This follows easily from the inequality
$$|c(x,y,z)-c(x_0,y,z)|\leq
C\frac{R}{|x-y|\,|x-z|},
$$
for $x,x_0,y,z$ such that
$|x-x_0|\leq R$ and $|x-y|,\,|x-z|\geq 2R$ (see Lemma 2.4 of
\cite{Tolsa1}, for example) and some standard estimates.
Therefore, if we suppose $M\sigma(x) \leq 1/(100C_3)$, then we
obtain
$$c_{\sigma|\C\setminus 2B}(x) \geq  c_{\sigma|\C\setminus 2B}(x_0)
- \frac1{100} \geq \frac1{10}.$$
So we have proved that if $M\sigma(x) \leq \min(1/1000,\,1/(6C_2)^{1/2},\,
1/(100C_3))$, then $c_\sigma(x) > 1/10$. This implies that in any case
we have $U_\sigma(x) \geq \alpha$, for  some $\alpha>0$.
\end{proof}


\section{The First Main lemma}  \label{sec4}

The proof of Theorem \ref{semiadd} uses an induction argument
on scales, analogous to the one in \cite{MTV}.
Indeed, if $Q$ is a sufficiently small square, then $E\cap Q$ either coincides
with a segment or it is void, and so
\begin{equation} \label{eqinduc}
\gamma_+(E\cap Q) \approx \gamma(E\cap Q).
\end{equation}
Roughly speaking, the idea consists of proving \rf{eqinduc}
for squares\footnote{Actually, in the induction argument we will use
rectangles instead of squares.}
$Q$ of any size, by induction.
To prove that \rf{eqinduc} holds for some fixed square $Q_0$,
we will take into account that \rf{eqinduc} holds for squares with side length
$\leq \ell(Q_0)/5$.

Our next objective consists of proving the following result.

\begin{lemma}[\bf First Main Lemma] \label{main1}
Suppose that $\gamma_+(E) \leq C_4\diam(E)$, with $C_4>0$ small enough.
Then there exists a compact set
$F=\bigcup_{i\in I}Q_i$, with $\sum_{i\in I}\chi_{10Q_i}\leq C$, such that
\begin{itemize}
\item[(a)] $E\subset F$ and $\gamma_+(F) \leq C\gamma_+(E)$,
\item[(b)] $\sum_{i\in I} \gamma_+(E\cap 2Q_i) \leq C\gamma_+(E)$,
\item[(c)] $\diam(Q_i)\leq \frac{1}{10}\diam(E)$ for every $i\in I$.
\end{itemize}
Let $A\geq1$ be some fixed constant and $\DD$ any fixed dyadic
lattice.
Suppose that $\gamma(E\cap 2Q_i)\leq A\gamma_+(E\cap 2Q_i)$ for all $i\in I$.
If $\gamma(E) \geq A\gamma_+(E)$, then there exist a positive Radon
measure $\mu$ and a complex Radon measure $\nu$, both supported on $F$,
and a subset $H_\DD\subset F$, such that
\begin{itemize}
\item[(d)] $C_a^{-1}\gamma(E) \leq \mu(F) \leq C_a\,\gamma(E)$,
\item[(e)] $d\nu = b\,d\mu$, with $\|b\|_{L^\infty(\mu)}\leq C_b$,
\item[(f)] $|\nu(F)| = \gamma(E),$
\item[(g)] $\int_{F\setminus H_\DD} \CC_*\nu\,d\mu \leq C_c\mu(F)$,
\item[(h)] If $\mu(B(x,r)) > C_0 r$ (for some big constant $C_0$), then
$B(x,r)\subset H_\DD$. In particular, $\mu(B(x,r))\leq C_0r$ for all
$x\in F\setminus H_\DD$, $r>0$.
\item[(i)] $H_\DD = \bigcup_{k\in I_H} R_k$, where $R_k$, $k\in I_H$, are
disjoint dyadic squares from the lattice $\DD$, with $\sum_{k\in I_H}
\ell(R_k)\leq \ve \mu(F)$, for $0<\ve<1/10$ arbitrarily
small (choosing $C_0$ big enough).
\item[(j)] $|\nu(H_\DD)| \leq \ve |\nu(F)|$.
\item[(k)] $\mu(H_\DD) \leq \delta \mu(F)$, with $\delta=\delta(\ve)<1$.
\end{itemize}
The constants $C_4$, $C$, $C_a$, $C_b$, $C_c$, $C_0$, $\ve$, $\delta$
do not depend on $A$. They are absolute constants.
\end{lemma}

Let us remark that the construction of the set $H_\DD$ depends on the
chosen dyadic lattice $\DD$. On the other hand, the construction of
$F$, $\mu$, $\nu$ and $b$ is independent of $\DD$.

We also insist on the fact that all the constants different from $A$
which appear in the lemma do not depend on $A$. This fact will be essential
for the proof of Theorem \ref{semiadd} in Section \ref{sec11}.
We have preferred to use the notation $C_a,C_b,C_c$ instead of $C_5,C_6,C_7$,
say, because these constants will play an important role in the proof of
Theorem \ref{semiadd}. Of course, the constant $C_b$ does not depend on $b$
(it is an absolute constant).

Remember that we said that we assumed the squares to be closed.
This is not the case for the squares of the dyadic squares of the lattice
$\DD$. We suppose that these squares are half open - half closed (i.e. of
the type $(a,b]\times(c,d]$).

For the reader's convenience, before going on we will make some
comments on the lemma. As we said in Section \ref{sketch}, the set $F$  has to
be understood as an approximation of $E$ at an intermediate scale.
The first part of the lemma, which deals with the construction of $F$ and the
properties (a)--(c), is proved in Section \ref{sec5}. The choice of the
squares $Q_i$ which satisfy (a) and (b) is one of the keys of the proof
Theorem \ref{semiadd}. Notice that (a) implies
that the squares $Q_i$ are not too big and (b) that they are not too
small. That is, they belong to some intermediate scale. The property (b)
will be essential for the proof of (d). On the other hand, the assertion (c)
will only be used in the induction argument, in Section \ref{sec11}.

The properties (d), (e), (f) and (g) are proved in Section \ref{sec6}.
These are the basic properties which must satisfy $\mu$ and $\nu$ in order to
apply a $T(b)$ theorem with absolute constants, as explained in Section
\ref{sketch}. To prove (d) we will need the assumptions in the paragraph
after (c) in the lemma. In (g) notice that instead of the
$L^\infty(\mu)$ or $BMO(\mu)$ norm of $\CC\nu$, we estimate the $L^1(\mu)$
norm of $\CC_*\nu$ out of the set $H_\DD$.

Roughly speaking, the {\em exceptional set} $H_\DD$ contains the part of $\mu$
without linear growth. The properties (h), (i), (j) and (k) describe $H_\DD$
and are proved in Section \ref{excep}. Observe that (i), (j) and (k) mean that
$H_\DD$ is a rather small set.


\section{Proof of {\rm (a)--(c)} in First Main Lemma}
\label{sec5}


\subsection{The construction of $F$ and the proof of (a)} \label{subsf}

Let $\sigma\in\Sigma(E)$ be a measure satisfying $\sigma(E)\approx
\gamma_+(E)$ and $U_\sigma(x)\geq\alpha>0$ for all $x\in E$ (recall Lemma
\ref{extremal}). Let $\lambda$ be some constant with $0<\lambda\leq
10^{-8}\alpha$ which will be fixed below.
Let $\Omega\subset\C$ be the {\bf open} set
$$\Omega:= \{x\in\C:\,U_\sigma(x)>\lambda\}.$$
Notice that $E\subset \Omega$, and by \cite[Theorem 3.1]{Tolsa2}
we have
\begin{equation}  \label{eqomega}
\gamma_+(\Omega) \leq C\lambda^{-1}\sigma(E) \leq C\lambda^{-1} \gamma_+(E).
\end{equation}

Let $\Omega=\bigcup_{i\in J} Q_i$ be a Whitney decomposition of $\Omega$,
where $\{Q_i\}_{i\in J}$ is the usual family of Whitney squares with
disjoint interiors satisfying
$20Q_i\subset\Omega$, $R Q_i\cap(\C\setminus\Omega)\neq \varnothing$ (where
$R$ is some fixed absolute constant), and $\sum_{i\in J}\chi_{10Q_i}\leq C$.

Let $\{Q_i\}_{i\in I}$, $I\subset J$,
be the subfamily of squares such that $2Q_i\cap E\neq
\varnothing$. We set
$$F:=\bigcup_{i\in I} Q_i.$$
Observe that the property (a) of First Main Lemma is a consequence
of \rf{eqomega} and the geometry of the Whitney decomposition.

To see that $F$ is compact it suffices to check that the family $\{Q_i\}_{i\in
I}$ is finite.
Notice that $E\subset \bigcup_{i\in J}(1.1\interior{Q}_i)$. Since
$E$ is compact, there exists a finite covering
$$E\subset \bigcup_{1\leq k\leq n}(1.1\interior{Q}_{i_k}).$$
Each square $2Q_i$, $i\in I$, intersects some square $1.1Q_{i_k}$,
$k=1,\ldots,n$.
Because of the geometric properties of the Whitney decomposition, the number
of squares $2Q_i$ which intersect some fixed square $1.1Q_{i_k}$ is bounded
above by some constant $C_5$. Thus, the family $\{Q_i\}_{i\in
I}$ has at most $C_5n$ squares.


\subsection{Proof of (b)}

Let us see now that (b) holds if $\lambda$ has been chosen small enough.
We will show below that if $x\in E\cap 2Q_i$ for some $i\in I$, then
\begin{equation} \label{1qaz}
U_{\sigma|4Q_i}(x)>\alpha/4,
\end{equation}
assuming that $\lambda$ is small enough. This implies
$E\cap2Q_i\subset \{U_{\sigma|4Q_i}>\alpha/4\}$ and then,
by \cite[Theorem 3.1]{Tolsa2}, we have
$$\gamma_+(E\cap 2Q_i) \leq C \alpha^{-1}\sigma(4Q_i).$$
Using the finite overlap of the squares $4Q_i$, we deduce
$$\sum_{i\in I} \gamma_+(E\cap2Q_i) \leq C\alpha^{-1}\sum_{i\in I}
\sigma(4Q_i) \leq C\alpha^{-1} \sigma(E) \leq C\gamma_+(E).$$
Notice that in the last inequality, the constant $\alpha^{-1}$
has been absorbed by the constant $C$.

Now we have to show that \rf{1qaz} holds for $x\in E\cap 2Q_i$.
Let $z\in RQ_i\setminus\Omega$, so
that $\dist(z,Q_i)\approx\dist(\partial\Omega,Q_i) \approx \ell(Q_i)$
(where $\ell(Q_i)$
stands for the side length of $Q_i$). Since $M\sigma(z) \leq U_\sigma(z) \leq
\lambda$, we deduce that for any square $P$ with $\ell(P)\geq \ell(Q_i)/4$ and
$P\cap 2Q_i\neq\varnothing$, we have
\begin{equation} \label{qgran}
\frac{\sigma(P)}{\ell(P)} \leq C_6\lambda \leq 10^{-6}\alpha,
\end{equation}
where the constant $C_6$ depends on the Whitney decomposition (in particular,
on the constant $R$), and we assume that $\lambda$ has been chosen so small
that the last inequality holds.

Remember that $U_\sigma(x) >\alpha$. If
$M\sigma(x) >\alpha/2$, then
$$\frac{\sigma(Q)}{\ell(Q)} > \alpha/2$$
for some ``small'' square $Q$ contained in $4Q_i$, because the ``big'' squares
$P$ satisfy \rf{qgran}. So, $U_{\sigma|4Q_i}(x) > \alpha/2$.

Assume now that $M\sigma(x) \leq \alpha/2$. In this case, $c_\sigma(x)
> \alpha/2$. We decompose $c^2_\sigma(x) =: c^2(x,\sigma,\sigma)$ as follows:
\begin{eqnarray*}
c^2(x,\sigma,\sigma) & = & c^2(x,\sigma|4Q_i,\sigma|4Q_i) +
2c^2(x,\sigma|4Q_i,\sigma|\C\setminus4Q_i) \\
&& \mbox{} +  c^2(x,\sigma|\C\setminus4Q_i,
\sigma|\C\setminus4Q_i).
\end{eqnarray*}
We want to see that
\begin{equation} \label{me3}
c_{\sigma|4Q_i}(x)>\alpha/4.
\end{equation}
So it is enough to show that the last two terms in the equation above are
sufficiently small.
First we deal with $c^2(x,\sigma|4Q_i,\sigma|\C\setminus4Q_i)$:
\begin{eqnarray} \label{me1}
c^2(x,\sigma|4Q_i,\sigma|\C\setminus4Q_i) & \leq & C\int_{y\in 4Q_i}
\int_{t\in\C\setminus4Q_i} \frac{1}{|t-x|^2}\,d\sigma(y)d\sigma(t)\nonumber\\
& = & C\mu(4Q_i) \int_{t\in\C\setminus4Q_i}
\frac{1}{|t-x|^2}\,d\sigma(t)\nonumber
\\ & \leq & C\mu(4Q_i) \frac{M\sigma(z)}{\ell(4Q_i)}
\leq C\,M\sigma(z)^2 \leq C\lambda^2.
\end{eqnarray}

For the term $c^2(x,\sigma|\C\setminus4Q_i,
\sigma|\C\setminus4Q_i)$ we write
\begin{eqnarray*}
c^2(x,\sigma|\C\setminus4Q_i,
\sigma|\C\setminus4Q_i) & = & c^2(x,\sigma|\C\setminus2RQ_i,
\sigma|\C\setminus2RQ_i) \\
&& \mbox{} + 2 c^2(x,\sigma|\C\setminus2RQ_i,
\sigma|2RQ_i\setminus4Q_i) \\
&&\mbox{} + c^2(x,\sigma|2RQ_i\setminus4Q_i,\sigma|2RQ_i\setminus4Q_i).
\end{eqnarray*}
Arguing as in \rf{me1}, it easily follows that the last two terms
are bounded above by $CM\sigma(z)^2\leq C\lambda^2$ again. So we get,
\begin{equation} \label{me2}
c^2_{\sigma|4Q_i}(x) \geq c^2_\sigma(x) - c^2_{\sigma|\C\setminus2RQ_i}(x)
- C\lambda^2.
\end{equation}
We are left with the term $c^2_{\sigma|\C\setminus2RQ_i}(x)$.
Since $x,z\in RQ_i$, it is not difficult to check
that
$$|c_{\sigma|\C\setminus 2RQ_i}(x) - c_{\sigma|\C\setminus 2RQ_i}(z)|
\leq C M\sigma(z) \leq C_6\lambda$$ [this is proved like
\rf{xx1}]. Taking into account that $c_\sigma(z) \leq \lambda$, we
get
$$c_{\sigma|\C\setminus 2RQ_i}(x) \leq (1+C_6)\lambda.$$
Thus, by \rf{me2}, we obtain
$$c^2_{\sigma|4Q_i}(x) \geq \frac{\alpha^2}{4} - C\lambda^2 \geq
\frac{\alpha^2}{16},$$
if $\lambda$ is small enough. That is, we have proved \rf{me3}, and so
in this case \rf{1qaz} holds too.


\subsection{Proof of (c)}

Now we have to show that
\begin{equation}\label{tam1}
\diam(Q_i)\leq \frac{1}{10}\diam(E).
\end{equation}
This will allow the application of our induction argument.

It is immediate to check that
$$U_\sigma(x) \leq \frac{100\sigma(E)}{\dist(x,E)}$$
for all $x\not\in E$ (of course, $100$ is not the best constant here).
Thus, for $x\in \Omega\setminus E$ we have
$$\lambda<U_\sigma(x) \leq \frac{100\sigma(E)}{\dist(x,E)}.$$
Therefore,
$$\dist(x,E)\leq 100\lambda^{-1}\sigma(E) \leq C\lambda^{-1}\gamma_+(E)
\leq \frac{1}{100}\,\diam(E),$$
taking the constant $C_4$ in First Main Lemma small enough.
As a consequence, $\diam(\Omega)\leq \frac{11}{10}\diam(E)$.
Since $20Q_i\subset \Omega$ for each $i\in I$, we have
$$20\,\diam(Q_i) \leq \diam(\Omega) \leq \frac{11}{10}\diam(E),$$
which implies \rf{tam1}.


\section{Proof of {\rm (d)--(g)} in First Main Lemma}
\label{sec6}


\subsection{The construction of $\mu$ and $\nu$ and the proof of (d)--(f)}

It is easily seen that there exists a family of $\CC^\infty$
functions $\{g_i\}_{i\in
J}$ such that, for each $i\in J$, $\supp(g_i)\subset 2Q_i$, $0\leq g_i\leq 1$,
and $\|\nabla g_i\|_\infty \leq C/\ell(Q_i)$, so that $\sum_{i\in J} g_i = 1$
on $\Omega$. Notice that by the definition of $I$ in Subsection \ref{subsf},
we also have $\sum_{i\in I} g_i = 1$ on $E$.

Let $f(z)$ be the Ahlfors function of $E$, and consider the complex measure
$\nu_0$ such that $f(z) = \CC\nu_0(z)$ for
$z\not \in E$, with $|\nu_0(B(z,r))|\leq r$ for all $z\in\C,\,r>0$
(see \cite[Theorem 19.9]{Mattila}, for example). So we have
$$|\CC\nu_0(z)| \leq 1\quad \mbox{ for all $z\not\in E$},$$
and
$$\nu_0(E) = \gamma(E).$$
The measure $\nu$ will be a suitable modification of $\nu_0$.
As we explained in Section \ref{sketch}, the main
drawback of $\nu_0$ is that the only information that we have about
its variation $|\nu_0|$ is that $|\nu_0|= b_0\,d\HH^1_E$, with
$\|b_0\|_\infty\leq1$. This is a very bad estimate if we try to apply
some kind of $T(b)$ theorem in order to show that the Cauchy transform
is bounded (with absolute constants). The main advantage of $\nu$ over
$\nu_0$ is that we will have a much better estimate for the variation $|\nu|$.

First we define the measure $\mu$. For each $i\in I$, let $\Gamma_i$
be a circumference concentric with $Q_i$ and radius $\gamma(E\cap 2Q_i)/10$.
Observe that $\Gamma_i\subset \frac{1}{2}Q_i$ for each $i$. We set
$$\mu = \sum_{i\in I} \HH^1|\Gamma_i.$$
Let us define $\nu$ now:
$$\nu = \sum_{i\in I} \frac{1}{\HH^1(\Gamma_i)}
\int g_i\,d\nu_0 \cdot\HH^1|\Gamma_i.$$
Notice that $\supp(\nu)\subset\supp(\mu)\subset F$. Moreover, we have
$\nu(Q_i) = \int g_i\,d\nu_0$, and since $\sum_{i\in I} g_i=1$ on $E$, we
also have $\nu(F) = \sum_{i\in I} \nu(Q_i) = \nu_0(E) = \gamma(E)$ (which
yields (f)).

We have $d\nu = b\,d\mu$, with $b =
\frac{\int g_i\,d\nu_0}{\HH^1(\Gamma_i)}$ on $\Gamma_i$. To estimate
$\|b\|_{L^\infty(\mu)}$, notice that
\begin{equation} \label{spli}
|\CC(g_i\nu_0)(z)|\leq C\quad \mbox{ for all $z\not\in E\cap 2Q_i$}.
\end{equation}
This follows easily from the formula
\begin{equation} \label{vit**}
\CC(g_i\nu_0)(\xi) = \CC\nu_0(\xi)\,g_i(\xi) + \frac{1}{\pi} \int
\frac{\CC\nu_0(z)}{z-\xi}\,\bar{\partial}g_i(z)\, d\LL^2(z),
\end{equation}
where $\LL^2$ stands for the planar Lebesgue measure on $\C$. Let
us remark that this identity is used often to split singularities
in Vitushkin's way. Inequality \rf{spli} implies that
\begin{equation} \label{xu}
\left|\int g_i\,d\nu_0\right| = \left|\left(\CC(g_i\,\nu_0)\right)'(\infty)
\right| \leq C \gamma(E\cap 2Q_i) = C\HH^1(\Gamma_i).
\end{equation}
As a consequence, $\|b\|_{L^\infty(\mu)}\leq C$, and (e) is proved.

It remains to check that (d) also holds.
Using \rf{xu}, the assumption $\gamma(E\cap 2Q_i)\leq A\gamma_+(E\cap2Q_i)$,
(b), and the hypothesis $A\gamma_+(E)\leq \gamma(E)$, we obtain the following
inequalities:
\begin{eqnarray*}
\gamma(E) = |\nu_0(E)| = \left|\sum_{i\in I} \int g_i\,d\nu_0\right|
& \leq & \sum_{i\in I} \left| \int g_i\,d\nu_0\right| \\
& \leq & C\sum_{i\in I} \gamma(E\cap 2Q_i) = C\mu(F) \\
& \leq & C\,A\sum_{i\in I} \gamma_+(E\cap 2Q_i)\\
& \leq & C\,A\gamma_+(E) \leq C\gamma(E),
\end{eqnarray*}
which gives (d) (with constants independent of $A$).

Notice, by the way, that the preceding inequalities show that
$\gamma(E)\leq CA\gamma_+(E)$. This is not very useful for us, because
if we try to apply induction, at each step of the induction the
constant $A$ will be multiplied by the constant $C$.

On the other hand, since for each square $Q_i$ we have $\mu(F\cap Q_i) \leq
CA\gamma_+(E\cap 2Q_i) \leq CA\sigma(2Q_i)$, with $\sigma\in\Sigma(E)$, it
follows easily that
\begin{equation} \label{creixA}
\mu(B(x,r))\leq C\,A\,r \quad \mbox{ for all $x\in F,\,r>0$.}
\end{equation}
Unfortunately, for our purposes this is not enough. We would like to obtain
the same estimate without the constant $A$ on the right hand side, but we will
not be able. Instead, we will get it for all $x\in F$
out of a rather small exceptional set $H$.


\subsection{The exceptional set $H$}

Before constructing the dyadic exceptional set $H_\DD$, we will consider
a non dyadic version, which we will denote by $H$.

Let $C_0\geq100C_a$ be some fixed constant. Following \cite{NTV}, given $x\in
F$, $r>0$, we say that $B(x,r)$ is a {\bf non Ahlfors disk} if
$\mu(B(x,r))> C_0 r$. For a fixed $x\in F$,
if there exists some $r>0$ such that $B(x,r)$ is a non Ahlfors disk, then
we say that $x$ is a {\bf non Ahlfors point}.
For any $x\in F$, we denote
$$\RR(x) = \sup\{r>0:\,B(x,r)\mbox{ is a non Ahlfors disk}\}.$$
If $x\in F$ is an Ahlfors point, we set $\RR(x)=0$.
We say that $\RR(x)$ is the {\bf Ahlfors radius} of $x$.

Observe (d) implies that $\mu(F) \leq C_a \gamma(E) \leq C_a\gamma(F)
\leq C_a\diam(F)$. Therefore,
$$\mu(B(x,r)) \leq \mu(F) \leq C_a\diam(F)\leq 100C_a r$$
for $r\geq\diam(F)/100$. Thus $\RR(x)\leq\diam(F)/100$ for all $x\in F$.

We denote
$$H_0 = \bigcup_{x\in F,\, \RR(x)>0} B(x,\RR(x)).$$
By Vitali's $5r$-Covering Theorem there is a disjoint family
$\{B(x_h,\RR(x_h))\}_h$ such that
$H_0 \subset \bigcup_h B(x_h,5\RR(x_h))$. We denote
\begin{equation} \label{defh}
H = \bigcup_h B(x_h,5\RR(x_h)).
\end{equation}
Since $H_0\subset H$, all non Ahlfors disks are contained in $H$ and then,
\begin{equation}\label{proph1}
\dist(x,F\setminus H) \geq \RR(x)
\end{equation}
for all $x\in F$.

Since
$\mu(B(x_h,\RR(x_h))) \geq C_0\RR(x_h)$
for every $h$, we get
\begin{equation} \label{proph3}
\sum_h \RR(x_h) \leq \frac1{C_0} \sum_h \mu(B(x_h,\RR(x_h)))
\leq \frac1{C_0}\mu(F),
\end{equation}
with $\frac{1}{C_0}$ arbitrarily small (choosing $C_0$ big enough).


\subsection{Proof of (g)}

The dyadic exceptional set $H_\DD$ will be constructed in Section \ref{excep}.
We will have $H_\DD\supset H$ for any choice of $\DD$.
In this subsection we will show that
\begin{equation} \label{g'}
\int_{F\setminus H} \CC_*\nu\,d\mu \leq C_c\mu(F),
\end{equation}
which implies (g), provided $H_\DD\supset H$.

We will work with the regularized operators $\CT$ introduced at
the end of Section \ref{prelim}. Remember that $|\CC\nu_0(z)|\leq
1$ for all $z\not\in E$. Since $\LL^2(E)=0$, the same inequality
holds $\LL^2$-a.e. $z\in\C$. Thus, $|\CT\nu_0(z)|\leq1$ and
$\wt{\CC}_*\nu_0(z)\leq1$ for all $z\in\C,\,\ve>0$.

To estimate $\wt{\CC}_*\nu$, we will deal with the term
$\wt{\CC}_*(\nu-\nu_0)$. This will be the main point for the proof of \rf{g'}.

We denote $\nu_i:=\nu|Q_i$.

\begin{lemma} \label{mpoint}
For every $z\in\C\setminus4Q_i$, we have
\begin{equation} \label{zxc1}
\wt{\CC}_*(\nu_i-g_i\nu_0)(z)\leq
\frac{C\ell(Q_i)\,\mu(Q_i)}{\dist(z,2Q_i)^{2}}.
\end{equation}
\end{lemma}

Notice that $\int (d\nu_i-g_i\,d\nu_0)=0$. Then, using the smoothness of the
kernels of the operators $\CT$, $\ve>0$, by standard estimates it easily
follows
$$\wt{\CC}_*(\nu_i-g_i\nu_0)(z)\leq
\frac{C\ell(Q_i)\,(|\nu|(Q_i)+|\nu_0|(2Q_i))}{\dist(z,2Q_i)^2}.$$
This inequality is not useful for our purposes because to estimate
$|\nu_0|(2Q_i)$ we only can use $|\nu_0|(2Q_i)\leq \HH^1(E\cap
2Q_i)$. However, we don't have any control over $\HH^1(E\cap
2Q_i)$ (we only know that it is finite, by our assumptions on
$E$). The estimate \rf{zxc1} is much sharper.

\begin{proof}[Proof of the lemma]
We set $\alpha_i = \nu_i-g_i\nu_0$. To prove the lemma, we have to
show that
\begin{equation} \label{dife}
|\CT\alpha_i(z)|\leq
\frac{C\ell(Q_i)\,\mu(Q_i)}{\dist(z,2Q_i)^{2}}
\end{equation}
for all $\ve>0$.

Assume first $\ve\leq \dist(z,2Q_i)/2$. Since $|\CC\alpha_i(w)|\leq C$
for all $w\not\in\supp(\alpha_i)$ and $\alpha_i(\C)=0$, we have
$$|\CC\alpha_i(w)| \leq \frac{C\,\diam(\supp(\alpha_i))\,
\gamma(\supp(\alpha_i))}{\dist(w,\supp(\alpha_i))^2}$$
(see \cite[p.12-13]{Garnett}). Remember that
$$\supp(\alpha_i)\subset\Gamma_i \cup (E\cap 2Q_i)\subset 2Q_i.$$
Then we get
\begin{equation} \label{dif2}
|\CC\alpha_i(w)| \leq \frac{C\,\ell(Q_i)
\gamma(\Gamma_i\cup(E\cap 2Q_i))}{\dist(w,2Q_i)^2}.
\end{equation}
Moreover, we have
$$\gamma(\Gamma_i\cup(E\cap 2Q_i)) \leq C(\gamma(\Gamma_i) +
\gamma(E\cap 2Q_i)),$$
because semiadditivity holds for the special case
$\Gamma_i\cup(E\cap 2Q_i)$. This fact follows easily from Melnikov's result
about semiadditivity of analytic capacity for two compacts which are separated
by a circumference \cite{Melnikov1}.
Therefore, by the definition of $\Gamma_i$, we get
\begin{equation} \label{dif3}
\gamma(\Gamma_i\cup(E\cap 2Q_i)) \leq C\gamma(E\cap 2Q_i) = C\mu(Q_i).
\end{equation}
If $w\in B(z,\ve)$, then $\dist(w,2Q_i)\approx\dist(z,2Q_i)$. By
\rf{dif2} and \rf{dif3} we obtain
$$|\CC\alpha_i(w)| \leq \frac{C\,\ell(Q_i)
\mu(Q_i)}{\dist(z,2Q_i)^2}.$$ Making the convolution with
$\psi_\ve$, \rf{dife} follows for $\ve\leq \dist(z,2Q_i)/2$.

Suppose now that  $\ve>\dist(z,2Q_i)/2$. We denote $h=\psi_\ve
*\alpha_i$. Then we have
$$\CT\alpha_i = \psi_\ve * \frac{1}{z} * \alpha_i
= \CC(h\,d\LL^2).$$ Therefore,
\begin{equation} \label{dif5}
|\CT\alpha_i(z)| \leq \int \frac{|h(\xi)|}{|\xi-z|}\, d\LL^2(\xi)
\leq \|h\|_\infty \,[\LL^2(\supp(h))]^{1/2}.
\end{equation}
We have to estimate $\|h\|_\infty$ and $\LL^2(\supp(h))$. Observe
that, if we write $\ell_i=\ell(Q_i)$ and we denote the center of
$Q_i$ by $z_i$, we have
$$\supp(h)\subset \supp(\psi_\ve) + \supp(\alpha_i) \subset
B(0,\ve) + B(z_i,2\ell_i) = B(z_i,\ve + 2\ell_i).$$ Thus,
$\LL^2(\supp(h))\leq C\ve^2$, since $\ell_i\leq\ve$.

Let us deal with $\|h\|_\infty$ now. Let $\eta_i$ be a
$\CC^\infty$ function supported on $3Q_i$ which is identically $1$
on $2Q_i$ and such that $\|\nabla\eta_i\|_\infty\leq C/\ell_i$.
Taking into account that $\alpha_i(2Q_i)=0$, we have
\begin{multline*}
h(w)   =  \int \psi_\ve(\xi-w)\,d\alpha_i(\xi)
 =  \int \bigl(\psi_\ve(\xi-w)- \psi_\ve(z_i-w)\bigr)\,d\alpha_i(\xi) \\
 =  \frac{\ell_i}{\ve^3}\int \frac{\ve^3}{\ell_i}\bigl(
\psi_\ve(\xi-w)-
\psi_\ve(z_i-w)\bigr)\,\eta_i(\xi)\,d\alpha_i(\xi)
 =:  \frac{\ell_i}{\ve^3}\int \vphi_w(\xi)\eta_i(\xi)\,d\alpha_i(\xi).
\end{multline*}
We will show below that
\begin{equation}
\label{zzx16} \|\CC(\vphi_w\eta_i\,d\alpha_i)\|_{L^\infty(\C)}\leq
C.
\end{equation}
Let us assume this estimate for the moment. Since
$\CC(\vphi_w\eta_id\alpha_i)$ is analytic in
$\C\setminus\supp(\alpha_i)$, using \rf{dif3} we deduce
$$\left|\frac{\ell_i}{\ve^3}\int
\vphi_w(\xi)\eta_i(\xi)\,d\alpha_i(\xi)\right| \leq
\frac{C\ell_i}{\ve^3} \gamma(\Gamma_i\cup(E\cap 2Q_i)) \leq
\frac{C\ell_i\mu(Q_i)}{\ve^3}.
$$
Therefore,
$$\|h\|_\infty\leq \frac{C\ell_i\mu(Q_i)}{\ve^3}.$$

By \rf{dif5} and the  estimates on $\|h\|_\infty$ and
$\LL^2(\supp(h))$, we obtain
$$|\CT\alpha_i(z)|
\leq \frac{C\ell(Q_i)\mu(Q_i)}{\ve^2} \leq
\frac{C\ell(Q_i)\mu(Q_i)}{\dist(z,2Q_i)^{2}}.$$

It remains to prove \rf{zzx16}. Remember that $\CC\alpha_i$ is a
bounded function. By the identity \rf{vit**}, since
$\supp(\vphi_w\eta_i)\subset 3Q_i$, it is enough to show that
\begin{equation}
\label{zzx17} \|\vphi_w\eta_i\|_\infty \leq C
\end{equation}
and
\begin{equation}
\label{zzx18} \|\nabla(\vphi_w\eta_i)\|_\infty \leq
\frac{C}{\ell_i},
\end{equation}
For $\xi\in3Q_i$, we have
$$|\vphi_w(\xi)| = \frac{\ve^3}{\ell_i}\,
\bigl|\psi_\ve(\xi-w)-\psi_\ve(z_i-w)\bigr| \leq \ve^3\,
\|\nabla\psi_\ve\|_{\infty} \leq C,$$ which yields \rf{zzx17}.
Finally, \rf{zzx18} follows easily too:
$$\|\nabla(\vphi_w\eta_i)\|_\infty \leq \|\nabla
\vphi_w\|_\infty + \|\vphi_w\|_\infty \|\eta_i\|_\infty \leq
\frac{C}{\ell_i}.$$ We are done.
\end{proof}


Now we are ready to prove \rf{g'}. We write
\begin{eqnarray} \label{eqh1}
\int_{F\setminus H} \wt{\CC}_*\nu\,d\mu & \leq &
\int_{F\setminus H} \wt{\CC}_*\nu_0\,d\mu +
\int_{F\setminus H} \wt{\CC}_*(\nu-\nu_0)\,d\mu \nonumber\\
& \leq & C\mu(F\setminus H) + \sum_{i\in I}
\int_{F\setminus H} \wt{\CC}_*(\nu_i-g_i\nu_0)\,d\mu.
\end{eqnarray}
To estimate the last integral we use Lemma \ref{mpoint} and recall that
$\|\wt{\CC}_*(\nu_i-g_i\nu_0)\|_{L^\infty(\mu)}\leq C$:
\begin{equation} \label{eqh2}
\int_{F\setminus H} \wt{\CC}_*(\nu_i-g_i\nu_0)\,d\mu \leq
C\mu(4Q_i) + \int_{F\setminus (4Q_i\cup H)}
\frac{C\ell(Q_i)\,\mu(Q_i)}{\dist(z,2Q_i)^2}\, d\mu(z).
\end{equation}
Let $N\geq1$ be the least integer such that
$(4^{N+1}Q_i\setminus4^{N}Q_i)\setminus H\neq\varnothing$, and take some fixed
$z_0\in (4^{N+1}Q_i\setminus4^{N}Q_i)\setminus H$. We have
\begin{eqnarray*}
\int_{F\setminus(4Q_i\cup H)}\frac{1}{\dist(z,2Q_i)^2}\, d\mu(z)
 & = & \sum_{k=N}^\infty \int_{4^{k+1}Q_i\setminus4^kQ_i)\setminus H}\\
& \leq & C \sum_{k=N}^\infty \frac{\mu(4^{k+1}Q_i)}{\ell(4^{k+1}Q_i)^2}\\
& \leq & C \sum_{k=N}^\infty
\frac{\mu(B(z_0,2\ell(4^{k+1}Q_i)))}{\ell(4^{k+1}Q_i)^2}\\
& \leq & C \sum_{k=N}^\infty
\frac{C_0 \ell(4^{k+1}Q_i)}{\ell(4^{k+1}Q_i)^2}\\
& \leq & C\,C_0\frac{1}{\ell(4^NQ_i)} \leq
C\,C_0\frac{1}{\ell(Q_i)}.
\end{eqnarray*}
Notice that in the second inequality we have used that $z_0\in F\setminus H$,
and so $\mu(B(z_0,r))\leq C_0r$ for all $r$.
By \rf{eqh2}, we obtain
$$\int_{F\setminus H} \wt{\CC}_*(\nu_i-g_i\nu_0)\,d\mu \leq C\mu(4Q_i).$$
Thus, by the finite overlap of the squares $4Q_i$, $i\in I$, and \rf{eqh1}, we
get
\begin{equation} \label{eqh4}
\int_{F\setminus H} \wt{\CC}_*\nu\,d\mu \leq C\mu(F\setminus H) +
C\sum_{i\in I}\mu(4Q_i) \leq C\mu(F).
\end{equation}
Now, \rf{w1} relates $\CC_*\nu$ with $\wt{\CC}_*\nu$:
\begin{equation} \label{eqh5}
|\wt{\CC}_*\nu(z)-\CC_*\nu(z)| \leq C\,M\nu(z).
\end{equation}
By (e), if $z\in F\setminus H$, we have
$M\nu(z)\leq CM\mu(z) \leq C$. Thus \rf{eqh4} and \rf{eqh5} imply
$$\int_{F\setminus H} \CC_*\nu(z)\,d\mu(z) \leq C\mu(F).$$


\section{The exceptional set  $H_\DD$}
\label{excep}


\subsection{The construction of $H_\DD$ and the proof of (h)--(i)}

Remember that in \rf{defh} we defined
$H= \bigcup_h B(x_h,5\RR(x_h))$, where $\{B(x_h,\RR(x_h))\}_h$ is some
precise family of non Ahlfors disks. Consider the family of
dyadic squares
$\DD_H\subset \DD$ such that $R\in\DD_H$ if there exists some ball
$B(x_h,5\RR(x_h))$ satisfying
\begin{equation} \label{cond1}
B(x_h,5\RR(x_h))\cap R\neq\varnothing
\end{equation}
and
\begin{equation} \label{cond2}
10\RR(x_h)< \ell(R) \leq 20\RR(x_h).
\end{equation}

Notice that
\begin{equation} \label{cond3}
\bigcup_h B(x_h,5\RR(x_h)) \subset \bigcup_{R\in\DD_H} R.
\end{equation}
We take a subfamily of disjoint maximal squares $\{R_k\}_{k\in I_H}$
from $\DD_H$ such that
$$\bigcup_{R\in\DD_H} R = \bigcup_{k\in I_H} R_k,$$
and we define the {\bf dyadic exceptional set $H_\DD$} as
$$H_\DD = \bigcup_{k\in I_H} R_k.$$
Observe that \rf{cond3} implies $H\subset H_\DD$ and, since for
each ball $B(x_h,5\RR(x_h))$ there are at most four squares
$R\in\DD_H$ satisfying \rf{cond1} and \rf{cond2}, by \rf{proph3}, we obtain
$$\sum_{k\in I_H} \ell(R_k) \leq 80\sum_h \RR(x_h)
\leq \frac{80}{C_0}\,\mu(F) \leq \ve\mu(F),$$
assuming $C_0\geq 80\ve^{-1}$.


\subsection{Proof of (j)}

Remember that the squares from the lattice $\DD$ are half open - half closed.
The other squares, such as the squares $\{Q_i\}_{i\in I}$ which form $F$, are
supposed to be closed.
From the point of view of the measures $\mu$ and $\nu$, there is no difference
between both choices because $\mu(\partial Q) =
|\nu|(\partial Q)= 0$ for any square $Q$ (remember that
$\mu$ is supported on a finite union of circumferences).

We have
$$|\nu(H_\DD)| \leq \sum_{k\in I_H} |\nu(R_k)|,$$
because the squares $R_k$, $k\in I_H$, are pairwise disjoint.
On the other hand, from (i), we deduce
$$\sum_{k\in I_H} \ell(R_k) \leq \ve\mu(F) \leq C_a \ve |\nu(F))|,$$
with $\ve\to 0$ as $C_0\to\infty$.
So (j) follows from next lemma.

\begin{lemma} \label{clcl}
For all squares $R\subset\C$, we have
$$|\nu(R)| \leq C\ell(R),$$
where $C$ is some absolute constant.
\end{lemma}

To prove this result we will need a couple of technical lemmas.


\begin{lemma} \label{jkh}
Suppose that $\bar{C}_0$ is some big enough constant. Let $R\subset\C$ be a
square such that $\mu(R)> \bar{C}_0\ell(R)$.
If $Q_i$ is a Whitney square such that $2Q_i\cap R\neq\varnothing$, then
$\ell(Q_i)\leq \ell(R)/4$.
\end{lemma}

\begin{proof}  Let us see that if $\ell(Q_i)> \ell(R)/4$, then
$\mu(R)\leq \bar{C}_0\ell(R)$. We may assume $\mu(R) \geq100\ell(R)$.
Notice that $R\subset 9Q_i$ and, by Whitney's
construction, we have $\#\{j:Q_j\cap9Q_i\neq\varnothing\}\leq C$. Further,
$\ell(Q_j)\approx\ell(Q_i)$ for this type of squares.
Recall also that
the measure $\mu$ on each Whitney square $Q_j$ coincides
with $\HH^1|\Gamma_j$, where $\Gamma_j$ is a circumference contained in
$\frac12 Q_j$, and so $\mu(Q_j)\leq C\ell(Q_j)$ for each $j$.
Therefore,
$$\mu(R) \leq \sum_{j:Q_j\cap9Q_i\neq\varnothing} \mu(Q_j)
\leq C\sum_{j:Q_j\cap9Q_i\neq\varnothing} \ell(Q_j) \leq C\ell(Q_i).$$
So we only have to show that $\ell(Q_i)\leq C\ell(R)$.

Since $2Q_i\cap R\neq\varnothing$, there exists some Whitney square
$Q_j$ such that $Q_j\cap R\neq\varnothing$ and
$Q_j\cap2Q_i\neq\varnothing$. Since we are assuming $\mu(R) \geq100\ell(R)$,
we have $\ell(R)\geq\ve_0\ell(Q_j)$, where $\ve_0>0$ is some absolute
constant (for instance, $\ve_0=1/100$ would possibly work). Thus,
$\ell(Q_i)\approx\ell(Q_j)\leq C\ell(R)$.
\end{proof}

\begin{lemma} \label{jkl}
Let $R\subset\C$ be a square such that $\ell(Q_i)\leq \ell(R)/4$ for each
Whitney square $Q_i$ with $2Q_i\cap R\neq \varnothing$.
Let $L_R = \{h\in I: 2Q_h\cap\partial R\neq\varnothing\}$.
Then,
$$\sum_{h\in L_R}\ell(Q_h)\leq C\ell(R).$$
\end{lemma}

\begin{proof}
Let $L$ be one of the sides of $R$.
Let $\{Q_h\}_{h\in I_L}$ be the subfamily of Whitney squares such that
$2Q_h\cap L\neq\varnothing$. Since $\ell(Q_h)\leq \ell(R)/4$, we have
$\HH^1(4Q_h\cap L)\geq C^{-1}\ell(Q_h)$. Then, by the bounded overlap of the
squares $4Q_h$, we obtain
\begin{equation} \label{hy1}
\sum_{h\in I_L} \ell(Q_h) \leq C\sum_{h\in I_L} \HH^1(4Q_h\cap L) \leq
C\ell(R).
\end{equation}
\end{proof}


\begin{proof}[Proof of Lemma \ref{clcl}]
By Lemma \ref{jkh}, we may assume $\ell(Q_i)\leq \ell(R)/4$ if $2Q_i\cap R\neq
\varnothing$. Otherwise, $|\nu(R)| \leq C_b\mu(R) \leq C_b\bar{C}_0\ell(R)$.

From the fact that $\|\CC\nu_0\|_{L^\infty(\C)}\leq 1$, we deduce
$|\nu_0(R)| \leq C\ell(R)$. So we only have to estimate the difference
$|\nu(R) - \nu_0(R)|$.

Let $\{Q_i\}_{i\in I_R}$, $I_R\subset I$, be the subfamily of Whitney squares
such that $Q_i\cap R\neq\varnothing$, and let
$\{Q_i\}_{i\in J_R}$, $J_R\subset I$, be the Whitney squares such that
$Q_i\subset R$. We write
\begin{eqnarray*}
|\nu(R) - \nu_0(R)| & = &
\Bigl|\nu\Bigl(\bigcup_{i\in I_R} (Q_i\cap R)\Bigr) -
 \nu_0\Bigl(\bigcup_{i\in I_R}(Q_i\cap R)\Bigr)\Bigr| \\
& \leq &
 \Bigl|\nu\Bigl(\bigcup_{i\in I_R\setminus J_R} (Q_i\cap R)\Bigr) -
 \nu_0\Bigl(\bigcup_{i\in I_R\setminus J_R}(Q_i\cap R)\Bigr)\Bigr|\\
&& \mbox{} +
\Bigl|\nu\Bigl(\bigcup_{i\in J_R} Q_i\Bigr) -
 \nu_0\Bigl(\bigcup_{i\in J_R}Q_i\Bigr)\Bigr| \\
& = & A + B.
\end{eqnarray*}
First we deal with the term $A$. We have
\begin{eqnarray*}
A & = & \Bigl|\sum_{i\in I_R\setminus J_R}
\frac{\HH^1(\Gamma_i\cap R)}{\HH^1(\Gamma_i)}
\int g_i\,d\nu_0 - \sum_{i\in I_R\setminus J_R} \nu_0(Q_i\cap  R)\Bigr|\\
& \leq & \sum_{i\in I_R\setminus J_R} \Bigl| \int g_i\,d\nu_0 \Bigr|
+ \sum_{i\in I_R\setminus J_R} |\nu_0(Q_i\cap  R)|.
\end{eqnarray*}
Since $|\CC(g_i\nu_0)|\leq C$ and $|\CC\nu_0|\leq C$, we have
$$\Bigl| \int g_i\,d\nu_0 \Bigr|
+ |\nu_0(Q_i\cap  R)| \leq C\ell(Q_i) + C\HH^1(\partial(Q_i\cap R))
\leq C\ell(Q_i).$$
Thus, $A\leq C\sum_{i\in I_R\setminus J_R}\ell(Q_i)$.
Notice now that if $i\in I_R\setminus J_R$, then $Q_i\cap  R\neq\varnothing$
and
$Q_i\not\subset R$. Therefore, $Q_i\cap
\partial R\neq\varnothing$.
From Lemma \ref{jkl} we deduce $A\leq C\ell(R)$.

Let us turn our attention to $B$:
\begin{eqnarray*}
B & = & \biggl|\sum_{i\in J_R} \int g_i\,d\nu_0 - \int_{\bigcup_{i\in J_R} Q_i}
d\nu_0\biggr| \\
& = & \biggl|\int \Bigl(\sum_{i\in J_R} g_i - \chi_{\bigcup_{i\in J_R} Q_i}
\Bigr)\, d\nu_0\biggr| \\
& \leq & \sum_{j\in J_R} \biggl|
\int_{Q_j} \Bigl(\sum_{i\in J_R} g_i - 1 \Bigr)\, d\nu_0 \biggr|
+ \biggl|\int_{\C\setminus \bigcup_{j\in J_R}\!Q_j}\,\sum_{i\in J_R}
g_i\,d\nu_0\biggr|\\
& = & B_1 + B_2.
\end{eqnarray*}

We consider first $B_1$. If $\sum_{i\in J_R} g_i\not\equiv1$ on $Q_j$, we write
$j\in M_R$. In this case there exists some $h\in I\setminus J_R$ such that
$g_h\not\equiv0$ on $Q_j$. So $2Q_h\cap Q_j\neq\varnothing$, with $Q_h\not
\subset  R$. Thus, $2Q_h\cap\partial R\neq\varnothing$.
That is, $h\in L_R$.

For each $h\in L_R$ there are at most $C_8$
squares $Q_j$ such that $2Q_h\cap Q_j\neq\varnothing$. Moreover, for these
squares $Q_j$ we have $\ell(Q_j) \leq C\ell(Q_h)$. Then, by Lemma \ref{jkl},
we get
\begin{equation} \label{jkl2}
\sum_{j\in M_R}\ell(Q_j) \leq C\,C_8\sum_{h\in L_R}\ell(Q_h) \leq C\ell(R).
\end{equation}

Now we set
\begin{eqnarray*}
B_1 & = & \sum_{j\in M_R}
\biggl|\int_{Q_j} \Bigl(\sum_{i\in J_R} g_i - 1 \Bigr)\, d\nu_0 \biggr|\\
& \leq & \sum_{j\in M_R} \biggl(
\biggl|\int_{Q_j} \sum_{i\in J_R} g_i \,d\nu_0\biggr| + |\nu_0(Q_j)|\biggr).
\end{eqnarray*}
We have $|\nu_0(Q_j)|\leq C\ell(Q_j)$ and also
$$\biggl|\int_{Q_j} \sum_{i\in J_R} g_i \,d\nu_0\biggr| \leq
\sum_{i\in J_R}\biggl|\int_{Q_j} g_i \,d\nu_0\biggr| \leq C\ell(Q_j),$$
because $\#\{i\in J_R:\supp(g_i)\cap Q_j\neq \varnothing \}\leq C$ and
$|\CC(g_i\nu_0)|\leq C$ for each $i$. Thus, by \rf{jkl2}, we deduce
$$B_1 \leq C\ell(R).$$

Finally we have to estimate $B_2$. We have
$$B_2 \leq \sum_{i\in J_R} \biggl|\int_{\C\setminus \bigcup_{j\in J_R}Q_j}
g_i\,d\nu_0\biggr| = \sum_{i\in J_R} B_{2,i}.$$
Observe that if $B_{2,i}\neq0$, then $\supp(g_i)\cap \supp(\nu_0)
\cap \C\setminus \bigcup_{j\in J_R}Q_j\neq \varnothing$. As  a consequence,
$2Q_i\cap Q_h\neq\varnothing$ for some $h\in I\setminus J_R$. Since
$Q_i\subset R$ and $Q_h\not\subset R$, we deduce that
either $2Q_i\cap\partial R\neq\varnothing$ or
$Q_h\cap\partial R\neq\varnothing$. So either $i\in L_R$ or $h\in L_R$.
Taking into account that $\ell(Q_i)\approx\ell(Q_h)$, arguing as above we get
\begin{eqnarray*}
B_2 & \leq & C\sum_{i\in L_R} \ell(Q_i) + C\sum_{i\in J_R}\,\,\sum_{h\in
L_R:Q_h\cap2Q_i\neq\varnothing} \ell(Q_i) \\
& \leq & C\ell(R) + C\sum_{h\in L_R} \,\,\sum_{i\in
I:Q_h\cap2Q_i\neq\varnothing}\ell(Q_h) \leq C\ell(R).
\end{eqnarray*}
\end{proof}


\subsection{Proof of (k)} \label{proofk}

Let us see that (k) is a direct consequence of (j).
We have
$$|\nu(F\setminus H_\DD) \geq |\nu(F)| - |\nu(H_\DD)| \geq (1-\ve) |\nu(F)|.$$
By (d) and (f), we get
$$\mu(F) \leq C|\nu(F)| \leq \frac{C}{1-\ve}\,|\nu(F\setminus H_\DD)|.$$
Since $\|b\|_{L^\infty(\mu)} \leq C$, we have $|\nu(F\setminus H_\DD)| \leq
C\mu(F\setminus H_\DD)$. Thus,
$\mu(F) \leq \frac{C_9}{1-\ve}\,\mu(F\setminus H_\DD)$. That is,
$\mu(H_\DD)\leq
\delta\mu(F)$, with $\delta=1-\frac{1-\ve}{C_9}$.


\section{The Second Main Lemma}
\label{sec8}

The second part of the proof of Theorem \ref{semiadd} is based on the
$T(b)$ theorem of Nazarov, Treil and Volberg in \cite{NTV}.
The precise result that we will prove is the following.
We use the same notation of First Main Lemma \ref{main1}.

\begin{lemma}[\bf Second Main Lemma]\label{main2}
Assume that $\gamma_+(E) \leq C_4\diam(E)$,
$\gamma(E) \geq A\gamma_+(E)$, and
$\gamma(E\cap 2Q_i)\leq A\gamma_+(E\cap 2Q_i)$ for all $i\in I$.
Then there exists some subset $G\subset F$,
with $\mu(F)\leq C_{10}\mu(G)$, such that
$\mu(G\cap B(x,r))\leq C_0r$ for all $x\in G$, $r>0$, and
the Cauchy transform is bounded on $L^2(\mu|G)$ with
$\|\CC\|_{L^2(\mu|G),L^2(\mu|G)} \leq C_{11}$, where $C_{11}$ is some
absolute constant. The constants $C_0$, $C_4$, $C_{10}$, $C_{11}$ are
absolute constants, and do not depend on $A$.
\end{lemma}

We will prove this lemma in the next two sections. First,
in Section \ref{secT} we will introduce two exceptional sets
$S$ and $T_\DD$ such that
$\CC_*\nu$ will be uniformly bounded on $F\setminus S$ and $b$
will behave as a paraacretive function out of $T_\DD$.
In the same section we will introduce the ``suppressed'' operators of
Nazarov, Treil and Volberg.
In Section \ref{secsup}
we will describe which modifications
are required in the $T(b)$ theorem of \cite{NTV} to prove Second Main Lemma.


\section{The exceptional sets $S$ and $T_\DD$ and the suppressed operators
$\CC_\Theta$}
\label{secT}


\subsection{The exceptional set $S$}

The arguments in this subsection will be similar to the ones in
\cite{Tolsa1.5}.

We set
$$S_0 = \{x\in F: \CC_*\nu >\alpha\},$$
where $\alpha$ is some big constant which will be chosen below.
For the moment, let us say that $\alpha\gg C_0C_b,C_c$.
For $x\in S_0$, let
$$\ve(x) = \sup\{\ve: \ve>0,\,|\CE\nu(x)|>\alpha\}.$$
Otherwise, we set $\ve(x)=0$.
We define the {\bf exceptional set $S$} as
$$S = \bigcup_{x\in S_0} B(x,\ve(x)).$$

To show that $\mu(S\setminus H_\DD)$ is small we will use the following
result.

\begin{lemma} \label{tamt}
If $y\in S\setminus H_\DD$, then $\CC_*\nu(y) >\alpha - 8C_0 C_b$.
\end{lemma}

\begin{proof}
Observe that if $y\in S\setminus H_\DD$, then  $y\in B(x,\ve(x))$ for some
$x\in S_0$. Let $\ve_0(x)$ be such that $|\CC_{\ve_0(x)}\nu(x)|>\alpha$ and
$y\in B(x,\ve_0(x))$.
We will show that
\begin{equation}  \label{mz2}
|\CC_{\ve_0(x)} \nu(x) - \CC_{\ve_0(x)} \nu(y)| \leq 8C_0C_b,
\end{equation}
and we will be done.
We have
\begin{eqnarray} \label{mz3}
\lefteqn{|\CC_{\ve_0(x)} \nu(x) - \CC_{\ve_0(x)} \nu(y)| \leq } &&\nonumber\\
&&|\CC_{\ve_0(x)} (\nu|B(y,2\ve_0(x))) (x)|
+ |\CC_{\ve_0(x)} (\nu|B(y,2\ve_0(x))) (y)| \nonumber \\
&& \mbox{}+
|\CC_{\ve_0(x)} (\nu|\C\setminus B(y,2\ve_0(x))) (x) - \CC_{\ve_0(x)}
(\nu|\C\setminus B(y,2\ve_0(x))) (y)|.
\end{eqnarray}
Notice that the first two terms on the right hand side are bounded above by
$$\frac{|\nu|(B(y,2\ve_0(x)))}{\ve_0(x)} \leq
\frac{C_b\mu(B(y,2\ve_0(x)))}{\ve_0(x)} \leq 2C_0 C_b,$$
since $y\not\in H_\DD$. The last term on the right hand side of \rf{mz3} is
bounded above by
\begin{multline*}
\int_{\C\setminus B(y,2\ve_0(x))} \left|\frac{1}{z-x} - \frac{1}{z-y} \right|
\,d|\nu|(z)
= \int_{\C\setminus B(y,2\ve_0(x))} \frac{|x-y|}{|z-x||z-y|} \,d|\nu|(z)\\
\leq 2C_b\ve_0(x) \int_{\C\setminus B(y,2\ve_0(x))} \frac{1}{|z-y|^2}
\,d\mu(z),
\end{multline*}
where we have applied that $|x-y|\leq\ve_0(x)$ and $|z-x|\geq|z-y|/2$ in the
last inequality.
As $y\not\in H_\DD$, we have the following standard estimate:
\begin{eqnarray*}
\lefteqn{2C_b\ve_0(x) \int_{\C\setminus B(y,2\ve_0(x))} \frac{1}{|z-y|^2}
\,d\mu(z) } & & \\
& = & 2C_b\ve_0(x) \sum_{k=1}^{\infty}
\int_{2^k \ve_0(x)\leq|z-y|<2^{k+1} \ve_0(x)} \frac{1}{|z-y|^2} \,d\mu(z) \\
& \leq & 2C_b\ve_0(x) \sum_{k=1}^{\infty} \frac{\mu(B(y,2^{k+1} \ve_0(x)))}{
2^{2k} \ve_0(x)^2}  \\
& \leq & 4 C_0 C_b.
\end{eqnarray*}
So we get
$$|\CC_{\ve_0(x)} (\nu|\C\setminus B(y,2\ve_0(x))) (x) - \CC_{\ve_0(x)}
(\nu|\C\setminus B(y,2\ve_0(x))) (y)| \leq 4 C_0 C_b,$$
and \rf{mz2} holds.
\end{proof}

Choosing $\alpha$ big enough,
we will have $\alpha/2\geq 8C_0C_b$. Then, from the
preceding lemma, we deduce
\begin{equation} \label{wert1}
\mu(S\setminus H_\DD) \leq \frac{2}{\alpha} \int_{F\setminus H_\DD}
\CC_*\nu\,d\mu
\leq \frac{2C_c}{\alpha}\,\mu(F),
\end{equation}
which tends to $0$ as $\alpha\to\infty$.



\subsection{The suppressed operators $\CC_\Theta$} \label{supres}

Let $\Theta:\C\lra [0,\infty)$ be a
Lipschitz function with Lipschitz constant $\leq 1$.
We denote
$$K_\Theta (x,y) = \frac{\overline{x-y}}{|x-y|^2 + \Theta(x)
\Theta(y)}.$$
It is not difficult to check that $K_\Theta$ is Calder\'on-Zygmund kernel
\cite{NTV}. Indeed, we have
$$|K_\Theta (x,y)| \leq \frac{1}{|x-y|},$$
and
$$|\nabla_x K_\Theta(x,y)|+ |\nabla_y K_\Theta(x,y)| \leq \frac{8}{|x-y|^2}.$$
The following estimate also holds:
\begin{equation}\label{idntv}
|K_\Theta(x,y)| \leq \frac{1}{\max\{\Theta(x),\Theta(y)\}}.
\end{equation}
We set
$$\CC_{\Theta,\ve} \nu(x) = \int_{\C\setminus B(x,\ve)} K_\Theta(x,y)
\,d\nu(y).$$
The operator $\CC_{\Theta,\ve}$ is the ($\ve$-truncated)
{\bf $\Theta$-suppressed Cauchy transform}.
We also denote
$$\CC_{\Theta,*} \nu(x)= \sup_{\ve>0} \CC_{\Theta,\ve} \nu(x).$$

The following lemma is a variant of some estimates which appear in
\cite{NTV}. It is also very similar to \cite[Lemma 2.3]{Tolsa1.5}.

\begin{lemma} \label{nnou}
Let $x\in\C$ and $r_0\geq0$ be such that $\mu(B(x,r))\leq C_0\,r$ for $r\geq
r_0$ and
$|\CE\nu(x)|\leq \alpha$ for $\ve\geq r_0$. If
$\Theta(x)\geq \eta r_0$ for some $\eta>0$, then
\begin{equation} \label{aaa}
|\CC_{\Theta,\ve}\nu (x)|\leq C_\eta
\end{equation}
for all $\ve>0$, with $C_\eta$ depending only on $C_0$, $C_b$,
 $\alpha$ and $\eta$.
\end{lemma}

\begin{proof}
If $\ve\geq \eta^{-1}\Theta(x)$, then
\begin{eqnarray*}
|\CC_{\Theta,\ve} \nu(x) - \CE\nu(x)| & \leq &
\int_{|x-y|>\ve} \left| \frac{\overline{x-y}}{|x-y|^2+\Theta(x)\Theta(y)}
- \frac{\overline{x-y}}{|x-y|^2} \right| d|\nu|(y) \\
& \leq &
C_b\int_{|x-y|>\ve} \left| \frac{(\overline{x-y})\,\Theta(x)\,\Theta(y)}{
|x-y|^2\,(|x-y|^2+\Theta(x)\Theta(y))} \right| d\mu(y) \\
& \leq & C_b\int_{|x-y|>\ve} \frac{\Theta(x)\Theta(y)}{|x-y|^3} d\mu(y) \\
& \leq & C_b\int_{|x-y|>\ve} \frac{\Theta(x)(\Theta(x)+|x-y|)}{|x-y|^3}
d\mu(y)\\
& = &  C_b\Theta(x)^2 \int_{|x-y|>\ve} \frac{1}{|x-y|^3} d\mu(y) \\
& & \mbox{} + C_b\Theta(x) \int_{|x-y|>\ve} \frac{1}{|x-y|^2} d\mu(y).
\end{eqnarray*}
Since $\mu(B(x,r))\leq C_0\,r$ for $r\geq \ve$, it is easily checked that
$$\int_{|x-y|>\ve} \frac{1}{|x-y|^3} d\mu(y) \leq \frac{C}{\ve^2}$$ and
$$\int_{|x-y|>\ve} \frac{1}{|x-y|^2} d\mu(y) \leq \frac{C}{\ve},$$
where $C$ depends only on $C_0$.
Therefore
$$|\CC_{\Theta,\ve} \nu(x) - \CE\nu(x)| \leq  \frac{C\Theta(x)^2}{\ve^2} +
\frac{C\Theta(x)}{\ve} \leq 2C,$$
and so \rf{aaa} holds for $\ve\geq \eta^{-1}\Theta(x)$.

If $\ve < \eta^{-1}\Theta(x)$, then
\begin{eqnarray}
|\CC_{\Theta,\ve} \nu(x)| & \leq & C_b \int_{B(x,\eta^{-1}\Theta(x))}
|K_\Theta(x,y)| d\mu(y) \nonumber \\
&& \mbox{} +
\biggl|\int_{\C \setminus B(x,\eta^{-1}\Theta(x))}
K_\Theta(x,y) d\nu(y)\biggr|.   \label{klk}
\end{eqnarray}
To estimate the first integral on the right hand side we use the inequality
\rf{idntv} and the fact that
$$\mu(B(x,\eta^{-1}\Theta(x))) \leq C_0\eta^{-1}\Theta(x),$$
because $\eta^{-1}\Theta(x)\geq r_0$.
The second integral on the right hand side of \rf{klk} equals
$\CC_{\Theta,\eta^{-1}\Theta(x)}\nu(x)$. This term is bounded by some constant,
as shown in the preceding case.
\end{proof}

We denote $\Phi_{0,\DD}(x) = \dist(x,\C\setminus (H_\DD\cup S))$.
Obviously, $\Phi_{0,\DD}(x)=0$ if $x\not\in H_\DD\cup S$. Moreover,
$\Phi_{0,\DD}$ is a Lipschitz function with Lipschitz constant $1$.
On the other hand, $H_\DD\cup S$ contains all non Ahlfors disks and all
the balls $B(x,\ve(x))$, $x\in F$, and so
$$\Phi_{0,\DD}(x) \geq \max(\RR(x),\ve(x)).$$
From the construction of $S$ and the preceding lemma we deduce:

\begin{lemma}
Let $\Theta:\C\lra[0,+\infty)$ be a Lipschitz function with Lipschitz constant
1 such that $\Theta(x)\geq \eta\,\Phi_{0,\DD}(x)$ for all $x\in\C$
(where $\eta>0$ is some fixed constant). Then, $\CC_{\Theta,*}\nu(x) \leq
C_\eta$ for all $x\in F$.
\end{lemma}


\subsection{The exceptional set $T_\DD$}

Looking at conditions (d), (e) and (f) of First Main Lemma \ref{main1} one can
guess that
the function $b$ will behave as a paraacretive function on many squares from
the dyadic lattice $\DD$. We deal with this question in this subsection.

Let us define the exceptional set $T_\DD$. If a dyadic square
$R\in \DD$ satisfies
\begin{equation} \label{acreno}
\mu(R) \geq C_d |\nu(R)|,
\end{equation}
where $C_d$ is some big constant which will be
chosen below, we write $R\in \DD_T$.
Let $\{R_k\}_{k\in I_T} \subset \DD_T$
be the subfamily of disjoint maximal dyadic squares from $\DD_T$.
The {\bf exceptional set $T_\DD$} is
$$T_\DD = \bigcup_{k\in I_T} R_k.$$

We are going to show that
$\mu(F\setminus (H_\DD \cup S \cup T_\DD))$ is big. That is,
that it is comparable to $\mu(F)$.
We need to deal with the sets $H_\DD$ and $T_\DD$ simultaneously.
Both $H_\DD$ and $T_\DD$ have been defined as a union of dyadic
squares satisfying some precise conditions (remember the property (i)
for the dyadic squares $R_k$, $k\in I_H$).

Let $\{R_k\}_{k\in I_{HT}}$ be the subfamily of different maximal (and thus
{\em disjoint}) squares from
$$\{R_k\}_{k\in I_{H}}\cup\{R_k\}_{k\in I_{T}},$$ so that
$$H_\DD \cup T_\DD = \bigcup_{k\in I_{HT}} R_k.$$
From Lemma \ref{clcl}, \rf{acreno} and the property (i) in First Main Lemma
\ref{main1}, we get
\begin{eqnarray*}
|\nu(H_\DD \cup T_\DD)| & \leq &\sum_{k\in I_{HT}} |\nu(R_k)| \\
& \leq & \sum_{k\in I_{H}} |\nu(R_k)| + \sum_{k\in I_{T}} |\nu(R_k)| \\
& \leq & C\sum_{k\in I_{H}} \ell(R_k) + C_d^{-1} \sum_{k\in I_T} \mu(R_k) \\
& \leq & C_{12}\ve \mu(F) + C_d^{-1}\mu(F) \\
& \leq & C_a(C_{12}\ve + C_d^{-1})|\nu(F)|.
\end{eqnarray*}
So if we choose $\ve$ small enough and $C_d$ big enough, we obtain
$$|\nu(H_\DD \cup T_\DD)| \leq \frac{1}{2}\,|\nu(F)|.$$
Now we argue as in Subsection \ref{proofk} for proving (k). We have
$$|\nu(F\setminus (H_\DD\cup T_\DD))| \geq
|\nu(F)| - |\nu(H_\DD\cup T_\DD)| \geq \frac{1}{2}\, |\nu(F)|.$$
Therefore,
$$\mu(F)  \leq  C_a|\nu(F)| \leq 2C_a|\nu(F\setminus (H_\DD\cup T_\DD))|
\leq  2C_aC_b\mu(F\setminus (H_\DD\cup T_\DD)).$$
Thus,
$\mu(H_\DD\cup T_\DD) \leq \delta_1 \mu(F)$, with
$\delta_1 = 1-\frac1{2C_aC_b} <1$.

Let us remark that the estimates above are not valid if we
argue with the non dyadic exceptional set $H$. We would have troubles
for estimating $\nu(H\cup T_\DD)$, because $H$ and $T_\DD$ are not disjoint
in general. This is the main reason for considering the dyadic version
$H_\DD$ of the exceptional set $H$ in First Main Lemma.

Now we turn our attention to the set $S$.
In \rf{wert1} we obtained an estimate for $\mu(S\setminus H_\DD)$
in terms of the constant $\alpha$.
We set $\delta_2 = (\delta_1+1)/2$. Then we choose $\alpha$ such that
$$\mu(H_\DD\cup T_\DD) + \mu(S\setminus H_\DD) \leq \delta_2\mu(F).$$


\subsection{Summary}
In next lemma we summarize what we have shown in this section.

\begin{lemma}  \label{summ}
Assume that $\gamma_+(E) \leq C_4\diam(E)$,
$\gamma(E) \geq A\gamma_+(E)$, and
$\gamma(E\cap Q)\leq A\gamma_+(E\cap Q)$ for all squares $Q$ with
$\diam(Q)\leq\diam(E)/5$.
Let $\DD$ be any fixed dyadic lattice.
There are subsets $H_\DD,\,S,\,T_\DD\subset F$
(with $H_\DD$ and $T_\DD$ depending
on $\DD$) such that
\begin{itemize}
\item[(a)] $\mu(H_\DD\cup S\cup T_\DD) \leq \delta_2\mu(F)$ for some absolute
constant $\delta_2<1$.

\item[(b)] All non Ahlfors disks (with respect to some constant $C_0$ big
enough) are contained in $H_\DD$.

\item[(c)] If $\Theta:\C\lra[0,+\infty)$ is any Lipschitz function with
Lipschitz constant 1 such that $\Theta(x)\geq \eta\,\dist(x,\C\setminus H_\DD
\cup S)$, for all $x\in\C$ (where $\eta>0$
is some fixed constant), then $\CC_{\Theta,*}\nu(x) \leq C_\eta$ for all
$x\in F$.

\item[(d)] All dyadic squares $R\in\DD$ such that $R\not\subset T_\DD$ satisfy
$\mu(R) < C_d |\nu(R)|$.
\end{itemize}
\end{lemma}


\section{The proof of Second Main Lemma}  \label{secsup}


Throughout all this section we will assume that all the hypotheses in
Second Main Lemma \ref{main2} hold.

\subsection{Random dyadic lattices}

We are going to introduce random dyadic lattices.
We follow the construction of \cite{NTV}.

Suppose that $F \subset B(0,2^{N-3})$, where $N$ is a big enough integer.
Consider the random square $Q^0(w) = w +\left[-2^N,\,2^N\right)^2$,
with $w\in \left[-2^{N-1},\,2^{N-1}\right)^2 =:\Omega$. We take $Q^0(w)$ as the
starting square of the dyadic lattice $\DD(w)$. Observe that $F\subset Q^0(w)$
for all $w\in\Omega$.
Only the dyadic squares which are
contained in $Q^0(w)$ will play some role in the arguments below. For the
moment, we don't worry about the other squares.

We take a uniform probability on $\Omega$.
So we let the probability measure $P$ be the normalized Lebesgue measure
on the square $\Omega$.

A square $Q\in\DD\equiv\DD(w)$ contained in $Q^0$ is called {\bf terminal} if
$Q\subset H_\DD\cup T_\DD$. Otherwise, it is called {\bf transit}.
The set of terminal squares is denoted by $\DD^{term}$, and the set of transit
squares by $\DD^{tr}$. It is easy to check that $Q^0$ is always transit.


\subsection{The dyadic martingale decomposition}

For $f\in L^1_{loc}(\mu)$ (we assume always $f$ real, for simplicity)
and any square $Q$ with $\mu(Q)\neq0$, we set
$$\langle f\rangle_Q = \frac1{\mu(Q)} \int_Q f\,d\mu.$$

We define the operator $\Xi$ as
$$\Xi f = \frac{\langle f\rangle_{Q^0}}{\langle b\rangle_{Q^0}}\,b,$$
where $b$ is the complex function that we have constructed in Main Lemma
\ref{main1}.
It follows easily that $\Xi f\in L^2(\mu)$ if $f\in L^2(\mu)$, and $\Xi^2
= \Xi$. Moreover, the definition of $\Xi$ does not depend on the choice
of the lattice $\DD$. The adjoint of $\Xi$ is
$$\Xi^*f = \frac{\langle fb\rangle_{Q^0}}{\langle b\rangle_{Q^0}}.$$

Let $Q\in \DD$ be some fixed dyadic square. The set of the four children of
$Q$ is denoted as $\CH(Q)$. In this subsection we will also write
$\CH(Q)= \{Q_j: j=1,2,3,4\}$.

For any square $Q\in\DD^{tr}$ and any $f\in L^1_{loc}(\mu)$,
we define the function $\D_Q f$ as follows:
$$\D_Q f = \left\{ \begin{array}{cl}
0 & \mbox{in $\C\setminus Q$,}\\
\ds \biggl(\frac{\langle f\rangle_{Q_j}}{\langle b\rangle_{Q_j}}
-\frac{\langle f\rangle_Q}{\langle b\rangle_{Q}}\biggr)b
& \mbox{in $Q_j$ if $Q_j\in\CH(Q)\cap\DD^{tr}$,} \\
\ds f - \frac{\langle f\rangle_Q}{\langle b\rangle_{Q}}\,b
& \mbox{in $Q_j$ if $Q_j\in\CH(Q)\cap\DD^{term}$.}
\end{array} \right.
$$

The operators $\D_Q$ satisfy the following properties.

\begin{lemma} \label{propde}
For all $f\in L^2(\mu)$ and all $Q\in\DD^{tr}$,
\begin{itemize}
\item[(a)] $\D_Q f\in L^2(\mu)$,
\item[(b)] $\int \D_Q f\,d\mu =0$,
\item[(c)] $\D_Q$ is a projection, i.e. $\D_Q^2=\D_Q$,
\item[(d)] $\D_Q\,\Xi = \Xi\,\D_Q =0$,
\item[(e)] If $R\in\DD^{tr}$ and $R\neq Q$, then $\D_Q\D_R=0$.
\item[(f)] The adjoint of $\D_Q$ is
$$\D_Q^* f = \left\{ \begin{array}{cl}
0 & \mbox{in $\C\setminus Q$,}\\
\ds \frac{\langle fb\,\rangle_{Q_j}}{\langle b\rangle_{Q_j}}
-\frac{\langle fb\rangle_Q}{\langle b\rangle_{Q}}
& \mbox{in $Q_j$ if $Q_j\in\CH(Q)\cap\DD^{tr}$,} \\
\ds f - \frac{\langle fb\rangle_Q}{\langle b\rangle_{Q}}
& \mbox{in $Q_j$ if $Q_j\in\CH(Q)\cap\DD^{term}$.}
\end{array} \right.
$$
\end{itemize}
\end{lemma}

The properties (a)--(e) are stated in \cite[Section XII]{NTV} and are
easily checked. The property (f) is also immediate (although it does not appear
in \cite{NTV}).

Now we have:

\begin{lemma} \label{decomp}
For any $f\in L^2(\mu)$, we have the decomposition
\begin{equation} \label{eqwa1}
f = \Xi f + \sum_{Q\in\DD^{tr}} \D_Q f,
\end{equation}
with the sum convergent in $L^2(\mu)$.
Moreover, there exists some absolute constant $C_{13}$ such that
\begin{equation} \label{eqwa2}
C_{13}^{-1}\|f\|_{L^2(\mu)}^2 \leq
\|\Xi f\|_{L^2(\mu)}^2 + \sum_{Q\in\DD^{tr}} \|\D_Q f\|_{L^2(\mu)}^2
\leq C_{13}\|f\|_{L^2(\mu)}^2.
\end{equation}
\end{lemma}

This lemma has been proved in \cite[Section XII]{NTV} under the assumption
that the paraacretivity constant $C_d$ (see \rf{acreno}) is sufficiently
close to $1$. The arguments in \cite{NTV} are still valid in our case for the
$L^2(\mu)$ decomposition of $f$
in \rf{eqwa1} and for the second inequality in \rf{eqwa2}.
However, they don't work for the first inequality in \rf{eqwa2}.
To prove it, we will use the Dyadic Carleson Imbedding Theorem:

\begin{theorem}
Let $\DD$ be some dyadic lattice and let $\{a_Q\}_{Q\in\DD}$ be a family of
non negative numbers. Suppose that for every square $R\in\DD$ we have
\begin{equation} \label{pack}
\sum_{Q\in\DD:Q\subset R} a_Q \leq C_{14}\mu(R).
\end{equation}
Then, for all $f\in L^2(\mu)$, we have
$$\sum_{Q\in\DD:\mu(Q)\neq0} a_Q\,|\langle f\rangle_Q|^2 \leq
4C_{14}\|f\|_{L^2(\mu)}^2.$$
\end{theorem}

See \cite[Section XII]{NTV}, for example, for the proof.

\begin{proof}[Proof of the first inequality in \rf{eqwa2}]
We will prove it by duality, like David in \cite{David3}. However we have to
modify the arguments because we cannot assume that $b^{-1}$ is a bounded
function (unlike in \cite{David3}),
since our function $b$ may vanish in sets of positive measure.

By \rf{eqwa1} and the fact that $\Xi$ and $\D_Q$ are projections, we have
$$f= \Xi f + \sum_{Q\in\DD^{tr}} \D_Q f =
\Xi^2 f + \sum_{Q\in\DD^{tr}} \D_Q^2 f.$$
Then we deduce
\begin{eqnarray} \label{kl2}
\int f^2 \,d\mu & = & \int \Bigl(\Xi^2 f + \sum_{Q\in\DD^{tr}} \D_Q^2 f\Bigr)
f\,d\mu \nonumber\\
& = & \int (\Xi f)(\Xi^* f)\,d\mu +  \sum_{Q\in\DD^{tr}}
\int (\D_Q f)(\D_Q^* f)\,d\mu \nonumber\\
& \leq &
\Bigl(\|\Xi f\|_{L^2(\mu)}^2 + \sum_{Q\in\DD^{tr}}
\|\D_Q f\|_{L^2(\mu)}^2\Bigr)^{1/2}\nonumber\\
&& \mbox{}\times
\Bigl(\|\Xi^* f\|_{L^2(\mu)}^2 + \sum_{Q\in\DD^{tr}}
\|\D_Q^* f\|_{L^2(\mu)}^2\Bigr)^{1/2}.
\end{eqnarray}
So if we show that
\begin{equation} \label{ws1}
\|\Xi^* f\|_{L^2(\mu)}^2 + \sum_{Q\in\DD^{tr}}
\|\D_Q^* f\|_{L^2(\mu)}^2 \leq C \|f\|_{L^2(\mu)}^2,
\end{equation}
we will be done. Notice, by the way, that the second inequality in \rf{eqwa2}
and \rf{kl2} imply
$$\|f\|_{L^2(\mu)}^2 \leq C \|\Xi^* f\|_{L^2(\mu)}^2 + C
\sum_{Q\in\DD^{tr}} \|\D_Q^* f\|_{L^2(\mu)}^2.$$

Let us see that \rf{ws1} holds. It is straightforward to check that
$$\|\Xi^*f\|_{L^2(\mu)} \leq C \|f\|_{L^2(\mu)}.$$
So we only have to estimate $\sum_{Q\in\DD^{tr}} \|\D_Q^* f\|_{L^2(\mu)}^2$.
To this end we need to introduce the operators $D_Q$. They are defined
as follows:
$$D_Q f = \left\{ \begin{array}{cl}
0 & \mbox{in $\C\setminus Q$,}\\
\langle f\rangle_{Q_j} - \langle f\rangle_Q
& \mbox{in $Q_j$.}
\end{array} \right.
$$
We also define $Ef = \langle f\rangle_{Q^0}$. Then it is well known that
\begin{equation} \label{ws44}
\|Ef\|_{L^2(\mu)}^2 + \sum_{Q\in\DD} \|D_Q f\|_{L^2(\mu)}^2 =
\|f\|_{L^2(\mu)}.
\end{equation}

If $Q_j\in\CH(Q)$ is a transit square, then we have (using (f) from Lemma
\ref{propde})
\begin{eqnarray*}
\D_Q^* f|_{Q_j} & = & \frac{\langle fb\rangle_{Q_j} - \langle fb\rangle_Q}{
\langle b\rangle_Q}  +
\langle fb\rangle_{Q_j} \biggl( \frac1{\langle b\rangle_{Q_j}} -
\frac1{\langle b\rangle_Q}\biggr) \\
& = & \frac1{\langle b\rangle_Q} D_Q(fb)|_{Q_j} -
\frac{\langle fb\rangle_{Q_j}}{\langle b\rangle_{Q_j} \langle b\rangle_Q}
D_Qb|_{Q_j}.
\end{eqnarray*}
Since $|\langle b\rangle_Q|,\, |\langle b\rangle_{Q_j}| \geq C_d^{-1}$,
we obtain
\begin{eqnarray} \label{ws45}
\lefteqn{\sum_{Q\in\DD^{tr}\,} \sum_{Q_j\in\CH(Q)\cap\DD^{tr}}
\|\D_Q^*f\|_{L^2(\mu|Q_j)}^2 \leq} &&\nonumber \\
&& C\sum_{Q\in\DD} \|D_Q(fb)\|_{L^2(\mu)}^2 +
C \sum_{Q\in\DD\,} \sum_{Q_j\in\CH(Q)}
\|\langle fb\rangle_{Q_j} D_Qb|_{Q_j}\|_{L^2(\mu|Q_j)}^2.
\end{eqnarray}
From \rf{ws44} we deduce
$$\sum_{Q\in\DD} \|D_Q(fb)\|_{L^2(\mu)}^2 \leq \|fb\|_{L^2(\mu)}^2 \leq
C\|f\|_{L^2(\mu)}^2.$$
Now observe that the last term in \rf{ws45} can be rewritten as
$$\sum_{Q\in\DD} |\langle fb\rangle_Q|^2
\|\chi_Q D_{\wh{Q}}b\|_{L^2(\mu)}^2 =:B,$$
where $\wh{Q}$ stands for the father of $Q$. To estimate this term we will
apply the Dyadic Carleson Imbedding Theorem. Let us check the numbers
$a_Q:= \|\chi_Q D_{\wh{Q}}b\|_{L^2(\mu)}^2$ satisfy the packing condition
\rf{pack}. Taking into account that $b$ is bounded and \rf{ws44}, for each
square $R\in\DD$ we have
\begin{eqnarray*}
\sum_{Q\subset R} \|\chi_Q D_{\wh{Q}}b\|_{L^2(\mu)}^2
& = & \|D_{\wh{R}}b\|_{L^2(\mu|R)}^2
+ \sum_{Q\subset R,\,Q\neq R} \|D_{\wh{Q}}b\|_{L^2(\mu)}^2 \\
& \leq & C\mu(R) + \sum_{Q\subset R} \|D_Q(b\chi_R)\|_{L^2(\mu)}^2
\leq C\mu(R).
\end{eqnarray*}
So \rf{pack} holds and then
$$B\leq C\|fb\|_{L^2(\mu)}^2 \leq C\|f\|_{L^2(\mu)}^2.$$

Now we have to deal with the terminal squares. If $Q\in\DD^{tr}$ and $Q_j
\in\DD^{term}$, then we have
$$\D_Q^* f|_{Q_j} = \biggl(f - \frac{\langle fb\rangle_{Q_j}}{\langle
b\rangle_Q}\biggr)
+ \biggl( \frac{\langle fb\rangle_{Q_j}}{\langle b\rangle_Q}
- \frac{\langle fb\rangle_{Q}}{\langle b\rangle_Q}\biggr).$$
Since $b$ is bounded and $|\langle b\rangle_Q|\geq C_d^{-1}$, we get
\begin{eqnarray*}
|\D_Q^* f|_{Q_j}| & \leq & C\bigl(|f| + \langle |f|\rangle_{Q_j}\bigr) +
C\bigl|\langle fb\rangle_{Q_j} - \langle fb\rangle_{Q}\bigr| \\
& = & C\bigl(|f| + \langle |f|\rangle_{Q_j}\bigr) + C
\bigl|D_Q(fb)|_{Q_j}\bigr|.
\end{eqnarray*}
Therefore,
\begin{eqnarray} \label{ws47}
\lefteqn{\sum_{Q\in\DD^{tr}\,} \sum_{Q_j\in\CH(Q)\cap\DD^{term}}
\|\D_Q^*f\|_{L^2(\mu|Q_j)}^2 \leq} &&\nonumber \\
&& C \sum_{Q\in\DD^{tr}\,} \sum_{Q_j\in\CH(Q)\cap\DD^{term}}
\int_{Q_j} |f|^2\,d\mu + C \sum_{Q\in\DD} \|D_Q (fb)\|_{L^2(\mu)}^2.
\end{eqnarray}
For the first sum on the right hand side above, notice
that the squares $Q_j\in\DD^{term}$ whose father is a transit square are
pairwise disjoint. For the last sum, we only have to use \rf{ws44}. Then
we obtain
$$\sum_{Q\in\DD^{tr}\,} \sum_{Q_j\in\CH(Q)\cap\DD^{term}}
\|\D_Q^*f\|_{L^2(\mu|Q_j)}^2 \leq C\|f\|_{L^2(\mu)}^2.$$
Since the left hand side of \rf{ws45} is also bounded above by
$C\|f\|_{L^2(\mu)}^2$, \rf{ws1} follows.
\end{proof}


\subsection{Good and bad squares}

Following \cite{NTV}, we say that a square $Q$ has {\bf $M$-negligible
boundary} if
$$\mu\{x\in\C:\dist(x,\partial Q)\leq r\} \leq Mr$$
for all $r\geq0$.

We now define bad squares as in \cite{NTV} too.
Let $\DD_1=\DD(w_1)$ and $\DD_2=\DD(w_2)$, with $w_1,w_2\in\Omega$, be two
dyadic lattices. We say that a transit square
$Q\in \DD_1^{tr}$ is bad (with respect to $\DD_2$) if either
\begin{itemize}
\item[(a)] there exists a square $R\in \DD_2$ such that
$\dist(Q,\partial R)\leq16\, \ell(Q)^{1/4}\ell(R)^{3/4}$ and $\ell(R) \geq
2^m\ell(Q)$ (where $m$ is some fixed positive integer), or

\item[(b)] there exists a square $R\in\DD_2$ such that
$R\subset (2^{m+2}+1)Q$, $\ell(R)\geq 2^{-m+1}\ell(Q)$, and $\partial R$ is not
$M$-negligible.
\end{itemize}
Of course, if $Q$ is not bad, then we say that it is good.

Let us remark that in the definition above we consider all the squares
of $\DD_2$, not only the squares contained in $Q^0(w_2)\in\DD_2$, which was
the case up to now.
On the other hand, observe that the definition depends on the constants $m$
and $M$. So strictly speaking, bad squares should be called
$(m,M)$-bad squares.

Bad squares don't appear very often in dyadic lattices. To be precise,
we have the following result.

\begin{lemma}[\cite{NTV}] \label{qbad}
Let $\ve_b>0$ be any fixed (small) number. Suppose that the constants $m$ and
$M$ are big enough (depending only on $\ve_b$). Let $\DD_1=\DD(w_1)$
be any fixed dyadic lattice.
For each fixed $Q\in\DD_1$, the probability that
it is bad with respect to a dyadic lattice $\DD_2=\DD(w_2)$, $w_2\in\Omega$,
is $\leq \ve_b$. That is,
$$P\{w_2:Q\in\DD_1 \mbox{ is bad with respect to $\DD(w_2)$}\} \leq\ve_b.$$
\end{lemma}

The notion of good and bad squares allows now to introduce the
definition of good
functions. Remember that given any fixed dyadic lattice $\DD_1=\DD(w_1)$,
every function $\vphi \in L^2(\mu)$ can be written as
$$\vphi = \Xi \vphi + \sum_{Q\in\DD_1^{tr}} \D_Q \vphi.$$
We say that $\vphi$ is {\bf $\DD_1$-good} with
respect to $\DD_2$ (or simply, good)
if $\D_Q\vphi=0$ for all bad squares $Q\in\DD_1^{tr}$
(with respect to $\DD_2$).


\subsection{Estimates on good functions}

We define the function $\Phi_\DD$ as
$$\Phi_\DD(x) = \dist(x,\C\setminus (H_\DD \cup S \cup T_\DD)).$$
Notice that $\Phi_\DD$ is a Lipschitz function with Lipschitz constant $1$
which equals zero in $\C\setminus (H_\DD \cup S \cup T_\DD)$.
Observe also that $\Phi_\DD\geq \Phi_{0,\DD}$ (this function was introduced
at the end of Subsection \ref{supres}).

Now we have the following result.

\begin{lemma} \label{estimgood}
Let $\DD_1=\DD(w_1)$ and $\DD_2=\DD(w_2)$, with $w_1,w_2\in\Omega$, be two
dyadic lattices.
Let $\Theta:\C\lra[0,+\infty)$ be a Lipschitz function with Lipschitz
constant 1 such that $\inf_{x\in\C}\Theta(x)>0$ and
$\Theta(x)\geq \eta\max(\Phi_{\DD_1}(x),\,\Phi_{\DD_2}(x))$
for all $x\in\C$ (where $\eta>0$ is some fixed constant).
If $\vphi$ is $\DD_1$-good with respect to $\DD_2$, and
$\psi$ is $\DD_2$-good with respect to $\DD_1$, then
$$|\langle \CC_\Theta \vphi,\,\psi\rangle| \leq C_{15} \|\vphi\|_{L^2(\mu)}
\|\psi\|_{L^2(\mu)},$$
where $C_{15}$ is some constant depending on $\eta$.
\end{lemma}

This lemma follows by
the same estimates and arguments of the corresponding result in \cite{NTV}.


\subsection{The probabilistic argument}

Following some ideas from \cite{NTV}, we are going to show that
the estimates for good functions from Lemma \ref{estimgood} imply
that there exists a set $G\subset \C\setminus H$ with
$\mu(G)\geq C^{-1}\mu(F)$ such that
the Cauchy transform is bounded on $L^2(\mu|G)$.
The probabilistic arguments of \cite[Section V]{NTV} don't work in our case
because we would need $\mu(H_\DD\cup S\cup T_\DD)$ to be very small
(choosing some adequate parameters), but we only have been able to show
that $\mu(H_\DD\cup S\cup T_\DD) \leq \delta_2\mu(F)$, for some fixed
$\delta_2<1$.
Nevertheless, the approach of \cite[Section XXIII]{NTV} doesn't need the
preceding assumption and is well suited for our situation.

Let us describe briefly the ideas from \cite[Section XXIII]{NTV}
that we need. We denote $W_\DD = H_\DD \cup S \cup T_\DD$, and we call it
the {\bf total exceptional set}.

Let $W_{\DD_1}$, $W_{\DD_2}$
be the total exceptional sets corresponding to two independent dyadic lattices
$\DD_1=\DD(w_1)$, $\DD_2=\DD(w_2)$. We have shown that
$$\mu(F\setminus W_{\DD(w)}) \geq
(1-\delta_2)\mu(F),$$ with $0<\delta_2<1$ for all $w\in\Omega$.
For each $x\in F$ we consider the probabilities
$$p_1(x) = P\{w\in\Omega:x\in F\setminus W_{\DD(w)}\},$$
and
$$p(x) = P\{(w_1,w_2)\in\Omega\times\Omega:x\in F\setminus (W_{\DD(w_1)}\cup
W_{\DD(w_2)})\}.$$
Since the sets $F\setminus W_{\DD(w_1)}$ and $F\setminus W_{\DD(w_2)}$ are
independent, we deduce $p(x) = p_1(x)^2$. Now we have
$$\int_F p_1(x)\,d\mu(x) = \EE\int\chi_{F\setminus W_{\DD(w)}}(x)\, d\mu(x)
=\EE\mu(F\setminus W_{\DD(w)}) \geq (1-\delta_2)\mu(F),$$
where $\EE$ denotes the mathematical expectation.
Let $G = \{x\in F: p_1(x)>(1-\delta_2)/2\}$, and $B=F\setminus G$. We have
\begin{eqnarray*}
\mu(B) & \leq & \frac{2}{1+\delta_2}\int_F (1-p_1(x))\,d\mu(x) \\
& = & \frac{2}{1+\delta_2} \biggl(\mu(F) - \int_F p_1(x)\,d\mu(x)\biggr)
\leq \frac{2\delta_2}{1+\delta_2}\,\mu(F).
\end{eqnarray*}
Thus,
$$\mu(G) \geq\frac{1-\delta_2}{1+\delta_2}\,\mu(F).$$
Observe that for every $x\in G$ we have $p(x) = p_1(x)^2 > (1-\delta_2)^2/4
=: \beta$. Now we define $\Phi_{(w_1,w_2)}(x) = \dist\bigl(x,
F\setminus(W_{\DD(w_1)}\cup W_{\DD(w_2)})\bigr)$. From the preceding
calculations, we deduce
$$\mu\{x\in F:p(x)>\beta\} \geq \mu(G) \geq
\frac{1-\delta_2}{1+\delta_2}\,\mu(F).$$
That is,
$$\mu\bigl\{x\in F: P\{(w_1,w_2):\Phi_{(w_1,w_2)}(x)=0\}>\beta\bigr\} \geq
\frac{1-\delta_2}{1+\delta_2}\,\mu(F).$$
Let us define
$$\Phi(x) = \inf_{B\subset\Omega\times\Omega,
P(B)=\beta\,\,}\sup_{\,(w_1,w_2)\in B}
\Phi_{(w_1,w_2)}(x).$$
Notice that $\Phi$ is a $1$-Lipschitz function such that $\Phi(x) =0$ for all
$x\in G$. Moreover, $\Phi(x)\geq \RR(x),\ve(x)$ for all $x\in F$, because
$\Phi_{(w_1,w_2)}(x)\geq \RR(x),\ve(x)$ for all $x\in
F,\,(w_1,w_2)\in\Omega\times\Omega$, since all non
Ahlfors disks are contained in $\HH_D$ for any choice of the lattice $\DD$,
and $S$ does not depend on $\DD$.

Finally, from Lemmas \ref{qbad} and \ref{estimgood}, and \cite[Main Lemma
(Section XXIII)]{NTV}, we deduce that $\CC_\Phi$ is bounded on
$L^2(\mu)$, and all the constants involved are absolute constants.
Since $\Phi(x)=0$ on $G$, {\bf the Cauchy transform is bounded
on $L^2(\mu|G)$}.
On the other hand, the fact that $\Phi(x)=0$ on $G$ also implies that
$\RR(x)=0$ on $G$, which is equivalent to say that $\mu(B(x,r))\leq C_0r$
for all $r>0$ if $x\in G$.

Now Second Main Lemma is proved.


\section{The proof of Theorem \ref{semiadd}}
\label{sec11}

From First Main Lemma and Second Main Lemma we get:

\begin{lemma} \label{lemaclau1}
There exists some absolute constant $B$ such that if $A\geq1$ is any fixed
constant and
\begin{itemize}
\item[(a)] $\gamma_+(E) \leq C_4\diam(E)$,
\item[(b)] $\gamma(E\cap Q)\leq A\gamma_+(E\cap Q)$ for all squares $Q$ with
$\diam(Q) \leq \diam(E)/5$,
\item[(c)] $\gamma(E) \geq A\gamma_+(E)$,
\end{itemize}
then $\gamma(E) \leq B\gamma_+(E)$.
\end{lemma}

\begin{proof} By First Main Lemma \ref{main1} and
Second Main Lemma \ref{main2}, there are sets $F,\,G$ and a measure $\mu$
supported on $F$ such that
\begin{enumerate}
\item[(1)] $E\subset F$ and $\gamma_+(E) \approx\gamma_+(F)$,
\item[(2)] $\mu(F) \approx\gamma(E)$,
\item[(3)] $G\subset F$ and $\mu(G)\geq C_{10}^{-1}\mu(F)$,
\item[(4)] $\mu(G\cap B(x,r))\leq C_0 r$ for all $x\in G$, $r>0$, and
$\|\CC\|_{L^2(\mu|G),L^2(\mu|G)} \leq C_{11}$.
\end{enumerate}
From (4) and (3), we get
$$\gamma_+(F) \geq C^{-1}\mu(G) \geq C^{-1}\mu(F).$$
Then, by (2), the preceding inequality, and (1),
$$\gamma(E) \leq C\mu(F) \leq C\gamma_+(F) \leq B\gamma_+(E).$$
\end{proof}

As a corollary we deduce:

\begin{lemma} \label{lemaclau2}
There exists some absolute constant $A_0$ such that if
$\gamma(E\cap Q)\leq A_0\gamma_+(E\cap Q)$ for all squares $Q$ with
$\diam(Q) \leq \diam(E)/5$, then
$\gamma(E) \leq A_0\gamma_+(E)$.
\end{lemma}

\begin{proof}
We take $A_0 = \max(1,C_4^{-1},B)$.
If $\gamma_+(E) > C_4 \diam(E)$, then
we get $\gamma_+(E) > C_4 \gamma(E)$ and we are done.
If $\gamma_+(E) \leq C_4 \diam(E)$, then
we also have $\gamma(E) \leq A_0\gamma_+(E)$. Otherwise,
we apply Lemma \ref{lemaclau1} and we deduce
$\gamma(E) \leq B\gamma_+(E) \leq A_0\gamma_+(E)$, which is a contradiction.
\end{proof}

Notice, by the way, that any constant $A_0\geq \max(1,C_4^{-1},B)$
works in the argument above. So Lemma \ref{lemaclau2}
holds for any constant $A_0$ sufficiently big.

Now we are ready to prove Theorem \ref{semiadd}.

\begin{proof}[\bf Proof of Theorem \ref{semiadd}]
Remember that we are assuming that $E$ is a finite union of disjoint compact
segments $L_j$. We set
$$d := \frac{1}{10}\,\min_{j\neq k}\dist(L_j,L_k).$$
We will prove by induction on $n$ that if $R$ is a closed rectangle with sides
parallel to the axes and diameter $\leq 4^n d$, $n\geq 0$, then
\begin{equation} \label{indu}
\gamma(R\cap E) \leq A_0 \gamma_+(R\cap E).
\end{equation}

Notice that if $\diam(R)\leq d$, then $R$ can intersect at most one segment
$L_j$. So either $R\cap E=\varnothing$ or $R\cap E$ coincides with a
segment, and in any case, \rf{indu} follows (assuming $A_0$ sufficiently big).

Let us see now that if
\rf{indu} holds for all rectangles $R$ with diameter $\leq 4^n d$,
then it also holds for a rectangle $R_0$ with diameter $\leq 4^{n+1} d$.
We only have to apply Lemma \ref{lemaclau2} to the set $R_0\cap E$, which
is itself a finite union of disjoint compact segments.
Indeed, take a square $Q$ with diameter $\leq \diam(R_0\cap E)/5$.
By the induction hypothesis we have
$$\gamma(Q\cap R_0\cap E) \leq A_0 \gamma_+(Q\cap R_0\cap E),$$
because $Q\cap R_0$ is a rectangle with diameter $\leq 4^nd$. Therefore,
$$\gamma(R_0\cap E) \leq A_0 \gamma_+(R_0\cap E)$$
by Lemma \ref{lemaclau2}.
\end{proof}


\end{document}